\documentclass[a4paper,12pt,twoside,notitlepage]{article}
\usepackage{a4,amsmath,amssymb}
\usepackage[noblocks]{authblk}
\usepackage{graphicx}
\def\q0{\theta}
\def\q{\vartheta}

\def\e0{\epsilon}

\def\f0{\phi}
\def\f{\varphi}

\def\r{\rho}

\def\R{{\mathbb R}}

\def\e{\varepsilon}

\newtheorem{thm}{Theorem}[section]
\newtheorem{cor}[thm]{Corollary}

\newtheorem{prop}[thm]{Proposition}
\newtheorem{defn}[thm]{Definition}
\newtheorem{rem}[thm]{Remark}
\begin{document}
\title{\sc The generalization of Sierpinski carpet and Sierpinski triangle in $n$-dimensional space\footnote{This work was supported by NSFC(Nos 11201056 and 11371080).}}
\author{Yun Yang, Yanhua Yu\thanks{Corresponding author. \\ {\it \hspace*{1em} Email addresses}: freeuse\_st@126.com (Yun Yang), yyh\_start@126.com(Yanhua Yu).}
\\{\it\small Department of Mathematics, Northeastern University,}\\{\it\small Shenyang, Liaoning, P. R. China,  110004}}
%\author{}
\markboth{}{}
\date{}
%\address{Huili Liu: Department of Mathematics, Northeastern University, Shenyang 110004, P. R. China}
%\email{liuhl@mail.neu.edu.cn, liu@math.tu-berlin.de}
%%%%%%%%%%%%%%%%%%%%%%%%%%%%%%%%
\maketitle
\numberwithin{equation}{section}
%%%%%%%%%%%%%%%%%%%%%%%%%%%%%%%%
\begin{abstract}
We obtain a nature generalization for an affine Sierpinski carpet and Sierpinski triangle to $n$-dimensional space, by using the generations and characterizations of affinely-equivalent Sierpinski carpet. Exactly, in this paper, a Menger sponge and Sierpinski simplex in $4$-dimensional space could be drawn out clearly under an affine transformation. Furthermore, the method could be used to a much broader class in fractals.
\medskip
\par
{\textbf{MSC 2000: }} 53A15, 81Q35.
\par
{\textbf{Key Words:}} Affine transformation, Sierpinski carpet, Menger sponge, affine invariants.
\end{abstract}
%%%%%%%%%%%%%%%%%%%%%%%%%%%%%%%%
% Begin the mail Text from now!
%%%%%%%%%%%%%%%%%%%%%%%%%%%%%%%%
\section{Introduction}
Discrete differential geometry studies discrete equivalents of the geometric notions and methods of classical differential
geometry, such as notions of curvature and integrability for polyhedral surfaces.  In this connection, discrete surfaces have been studied one after another with
strong ties to mathematics physics and great potential for computer analysis, architecture, numerics. Progress in this field is to a
large extent stimulated by its relevance for computer graphics and mathematical physics\cite{Bobenko-2,Bobenko-4}.
\par Recently, the expansion of computer graphics and applications in mathematical physics have given a great impulse to the issue of giving discrete equivalents of affine differential geometric objects\cite{Bobenko-1}. In \cite{Bobenko-3} a consistent definition of discrete affine spheres is proposed, both for definite and indefinite metrics, and in \cite{Matsuura} a similar construction is done in the context of improper affine spheres.
\par Following the ideas of Klein, presented in his famous lecture at Erlangen, several geometers in the early 20th century proposed the study
of curves and surfaces with respect to different transformation groups. In geometry, an affine transformation, affine map or an affinity is a function between affine spaces which preserves points, straight lines and planes. Also, sets of parallel lines remain parallel after an affine transformation. An affine transformation does not necessarily preserve angles between lines or distances between points, though it does preserve ratios of distances between points lying on a straight line. Examples of affine transformations include translation, scaling, homothety, similarity transformation, reflection, rotation, shear mapping, and compositions of them in any combination and sequence.
\par A centro-affine transformation is nothing but a general linear transformation $\R^n\ni x\mapsto Ax\in\R^n$, where $A\in GL(n,\R)$.
 In 1907 Tzitz$\mathrm{\acute{e}}$ica found that for a surface in Euclidean 3-space the property that the ratio of the Gauss curvature to the fourth power of
the distance of the tangent plane from the origin is constant is invariant under a centro-affine transformation. The surfaces with this property
turn out to be what are now called Tzitz$\mathrm{\acute{e}}$ica surfaces, or proper affine spheres with center at the origin. In centro-affine differential geometry, the theory of hypersurfaces has a long history. The notion of centro-affine minimal
hypersurfaces was introduced by Wang \cite{Wang} as extremals for the area integral of the centro-affine metric. See also \cite{Y-Y-L,Yu-Y-L} for the classification results about centro-affine translation surfaces and centro-affine ruled surfaces in $\R^3$.
\par Smooth geometric objects and their transformations should belong to the same geometry. In particular discretizations should be invariant with respect to the same transformation group as the smooth objects are(projective, affine, m\"obius etc). Deterministic and statistical fractals are an important tool for the investigation of physical phenomena. They were used by Mandelbrot to describe physical characteristics of things such as rivers, coastlines, bronchi and music. A vast literature on the theory and application of fractals has appeared. Generalization of the Sierpinski triangle to the Sierpinski tetrahedron in three dimensions has been described
in several places and generalizations to every dimension and every base for representing addresses have also been described. All of these generalizations are point
based in that the finite approximations can essentially be described as points in $n$-dimensional space. These representations are difficult to deal with visually in dimensions higher than three because of the large number of points and the difficulty in choosing a projection to preserve as much information and symmetry as possible\cite{Brisson}.

 \par In this paper we consider a nature generalization for an affine Sierpinski carpet to $n$-dimensional space, by using the generations and characterizations of affinely-equivalent Sierpinski carpet, which is organized as follows: Basic concepts of classical centro-affine geometry are presented in Section 2.
In Section 3 we define the discrete centro-affine hypersurface, and then obtain the structure equations, compatibility conditions and some centro-affine invariants.
In section 4, we study the generations and characters of Sierpinski carpet and Menger sponge, and find their regularities. Then we generalize the definition of  Menger sponge to $n$-dimensional space. Finally, we draw out the graphs of Menger sponge in $4$-dimensional space. In section 5, using the same method as in Section 4, we obtain the generalization of Sierpinski triangle to $n$-dimensional space.

\section{Affine mappings and transformation groups, basic notations}
If $X$ and $Y$ are affine spaces, then every affine transformation $f:X\rightarrow Y$  is of the form $\mathbf{x}\mapsto M\mathbf{x}+\mathbf{b}$ , where $M$ is a linear transformation on $X$ and  $b$ is a vector in $Y$. Unlike a purely linear transformation, an affine map need not preserve the zero point in a linear space. Thus, every linear transformation is affine, but not every affine transformation is linear.
\par For many purposes an affine space can be thought of as Euclidean space, though the concept of affine space is far more general (i.e., all Euclidean spaces are affine, but there are affine spaces that are non-Euclidean, for example, Minkowski space). In affine coordinates, which include Cartesian coordinates in Euclidean spaces, each output coordinate of an affine map is a linear function (in the sense of calculus) of all input coordinates. Another way to deal with affine transformations systematically is to select a point as the origin; then, any affine transformation is equivalent to a linear transformation (of position vectors) followed by a translation.
\par It is well known that the set of all automorphisms of a vector space $V$ of dimension $m$ forms a group. We use the following standard notations for this group and its subgroups(\cite{L-U-Z}):
$$GL(m,\R):=\{L:V\rightarrow V|L\quad isomorphism\};$$
$$SL(m,\R):=\{L\in GL(m,\R)|\det L=1\}.$$
Correspondingly, for an affine space $A, \dim A=m$, we have the following affine transformation groups.
\begin{flalign*}
% \nonumber to remove numbering (before each equation)
  &\mathcal{A}(m):=\{\alpha:A\rightarrow A|L_{\alpha}\ \mathrm{regular}\}\quad \mathrm{is\ the\ regular\ affine\ group}.  \\
  &\mathcal{S}(m):=\{\alpha\in \mathcal{A}|\det L_\alpha=1\}\quad \mathrm{is\ the\ unimodular(equiaffine)\ group}.  \\
  &\mathcal{Z}_p(m):=\{\alpha\in \mathcal{A}|\alpha(p)=p\}\quad \mathrm{is\ the\ centro-affine\ group\ with\ center}\ p\in A. \\
  &\tau(m):=\{\alpha:A\rightarrow A| \mathrm{there\ exists}\ b(\alpha)\in V, \mathrm{s.t.}\ \overrightarrow{p\alpha(p)}=b(\alpha), \forall p\in A\}\\
  &\qquad\qquad \mathrm{is\ the\ group\ of\ transformations\ on\ } A.
\end{flalign*}
Let $\mathcal{G}$ be one of the groups above and $S_1,S_2\subset A$ subsets. Then $S_1$ and $S_2$ are called equivalent modulo $\mathcal{G}$ if there exists an $\alpha\in\mathcal{G}$ such that $$S_2=\alpha S_1.$$
\par In centro-affine geometry we fix a point in $A$(the origin $O\in A$ without loss of generality) and consider the geometric properties in variant under the centro-affine group $\mathcal{Z}_p$. Thus the mapping $\pi_0:A\rightarrow V$ identifies $A$ with the vector space $V$ and $\mathcal{Z}_O$ with $GL(m, \R)$.
\section{Discrete centro-affine hypersurfaces.}
Here, we introduce discrete analogues of centro-affine hypersurfaces though a purely geometric manner. These constitute particular `discrete hypersurface' which are maps
\begin{equation}\label{dis-map-h}
  \mathbf{r}:\mathbb{Z}^n\rightarrow \mathbb{R}^{n+1},\qquad (k_1,k_2,\cdots,k_n)\mapsto \mathbf{r}(k_1,k_2,\cdots,k_n).
\end{equation}
In the following, we suppress the arguments of functions of $k_1,k_2,\cdots,k_n$, and denote increments of the discrete variables by subscripts, for example,
$$\mathbf{r}=\mathbf{r}(k_1,\cdots,k_n),\mathbf{r}_1=\mathbf{r}(k_1+1,\cdots,k_n),\mathbf{r}_n=\mathbf{r}(k_1,\cdots,k_n+1).$$
Moreover, decrements are indicated by overbars, that is,
$$\mathbf{r}_{\bar{1}}=\mathbf{r}(k_1-1,\cdots,k_n),\mathbf{r}_{\bar{n}}=\mathbf{r}(k_1,\cdots,k_n-1).$$
Our convention for the range of indices is the following
$$1\leq i,j,k,\cdots\leq n,$$
$$1\leq A,B,C,\cdots\leq n+1.$$
Now we will give a definition for the discrete centro-affine hypersurface. Especially, in the following, the point $\mathbf{r}(k_1,\cdots,k_n)$ indicates the terminal point of the vector $\mathbf{r}(k_1,\cdots,k_n)$ with its starting point at the origin $O$.
\begin{defn}\label{DCS-h}(Discrete centro-affine hypersurface) A $n$-dimensional lattice (net) in $n+1$-dimensional affine space
\begin{equation}
 \mathbf{r}:\mathbb{Z}^n\rightarrow \mathbb{R}^{n+1}
\end{equation}
is called a discrete centro-affine hypersurface if  the vector group $\mathbf{r},\mathbf{r}_1,\mathbf{r}_2,\cdots \mathbf{r}_n$ are linearly independent, and
the vector group $\mathbf{r},\mathbf{r}_{\bar{1}},\mathbf{r}_{\bar{2}},\cdots \mathbf{r}_{\bar{n}}$ are also linearly independent.
\end{defn}
 According to Definition \ref{DCS-h}, we know $\mathbf{r}, \mathbf{r}_1,\cdots, \mathbf{r}_n$ are linearly independent, so the following chain structures are obvious.
\begin{align}
% \nonumber to remove numbering (before each equation)
  \left\{\mathbf{r}_1,\mathbf{r}_{11},\cdots,\mathbf{r}_{1n}\right\}=&\left\{\mathbf{r},\mathbf{r}_1,\cdots,\mathbf{r}_n\right\}M1(k_1,\cdots,k_n), \label{Tran-1-h}\\
  \left\{\mathbf{r}_2,\mathbf{r}_{21},\cdots,\mathbf{r}_{2n}\right\}=&\left\{\mathbf{r},\mathbf{r}_1,\cdots,\mathbf{r}_n\right\}M2(k_1,\cdots,k_n),  \label{Tran-2-h}\\
  \cdots\quad & \quad \cdots\nonumber\\
  \left\{\mathbf{r}_n,\mathbf{r}_{n1},\cdots,\mathbf{r}_{nn}\right\}=&\left\{\mathbf{r},\mathbf{r}_1,\cdots,\mathbf{r}_n\right\}Mn(k_1,\cdots,k_n),  \label{Tran-3-h}
\end{align}
where $(n+1)\times(n+1)$ matrices $M1,M2,\cdots,Mn$ are non-degenerate. The first column of Matrix $Mk$ is
$$(0,\cdots,0,1,0,\cdots,0)^{\mathrm{Tran}},$$ and the element $1$ appears at position $k+1$. Especially, the $(i+1)$th column of Matrix $Mk$ is the same as the $(k+1)$ column of Matrix $Mi$, that is
$$(ei)_{A,k+1}=(ek)_{A,i+1},$$
where $(ei)_{A,k+1}$ represents the element of Matrix $Mi$ in row $A$ and column $k+1$.
Furthermore, Eqs. (\ref{Tran-1-h})-(\ref{Tran-3-h}) yields the compatibility conditions
\begin{equation}\label{com-con-h}
 (Mi)(Mj)_i = (Mj)(Mi)_j.
\end{equation}
Hence, we get the following structure equations and compatibility conditions from Eqs. (\ref{Tran-1-h})-(\ref{Tran-3-h}) and (\ref{com-con-h}).
\begin{prop}The structure equations of the discrete centro-affine hypersurface can be written as
\begin{eqnarray}
% \nonumber to remove numbering (before each equation)
   \mathbf{r}_{ii}-\mathbf{r}_i&=&(\sum_{k=1}^{n+1}(ei)_{k,i+1}-1)\mathbf{r}+(ei)_{2,i+1}(\mathbf{r}_1-\mathbf{r})+\cdots+((ei)_{i+1,i+1}-1)(\mathbf{r}_i-\mathbf{r})\nonumber\\
   & &\qquad +\cdots+(ei)_{n+1,i+1}(\mathbf{r}_n-\mathbf{r}), \label{Stru-1-h} \\
   \mathbf{r}_{ij}-\mathbf{r}&=&(\sum_{k=1}^{n+1}(ei)_{k,j+1}-1)\mathbf{r}+(ei)_{2,j+1}(\mathbf{r}_1-\mathbf{r})+\cdots+(ei)_{n+1,j+1}(\mathbf{r}_n-\mathbf{r}),   \label{Stru-2-h}
\end{eqnarray}
where $(ei)_{A,B}$ represents the element of Matrix $Mi$ in row $A$ and column $B$ and $i\neq j$.
Its compatibility conditions are
$$(Mi)(Mj)_i = (Mj)(Mi)_j.$$
\end{prop}
Indeed, from Eqs. (\ref{Tran-1-h})-(\ref{Tran-3-h}) we obtain
\begin{gather}
(ei)_{A,l+1}=\frac{[\mathbf{r},\mathbf{r}_1,\cdots,\mathbf{r}_{il},\cdots,\mathbf{r}_n]}{[\mathbf{r},\mathbf{r}_1,\cdots,\mathbf{r}_{A-1},\cdots,\mathbf{r}_n]},
\end{gather}
where $[\dots]$ is the standard determinant in $\R^{n+1}$. Hence, it is easy to conclude
\begin{cor} The discrete functions $(ei)_{A,B}$ are centro-affine invariant.
\end{cor}
Given the discrete functions $(ei)_{A,B}$, which satisfy the compatibility condition (\ref{com-con-h}), and then by using the initial $n+1$ points $\mathbf{r}$,$\mathbf{r}_1,\cdots,\mathbf{r}_n,$ all point groups in the discrete centro-affine surface can be generated by using the chain rule (\ref{Tran-1-h})-(\ref{Tran-3-h}). Therefore, we can state the following results.
\begin{cor}Two discrete centro-affine surface are centro-affinely equivalent if and only if the invariants $(ei)_{A,B}$ are same.
\end{cor}
\begin{cor}\label{self-sim}
A discrete centro-affine surface is self-similar if the centro-affine invariants $(ei)_{A,B}$ are constant.
\end{cor}
Especially, if a hypersurface lies on a hyperplane which does not contain the origion $O$,  from Eqs. (\ref{Stru-1-h})-(\ref{Stru-2-h}), it is easy to see the coefficients of position vector $\mathbf{r}$ should be zero.
And then, we can obtain that all coefficients are invariant under the translation transformation. Hence the following result are obvious.
\begin{prop}\label{prop-hypplan}The hypersurface lies in a hyperplane which does not contain the origion $O$ if and only if
 $$\sum_{k=1}^{n+1}(ei)_{k,i+1}=1, \quad \sum_{k=1}^{n+1}(ei)_{k,j+1}=1.$$ In this case, the functions $(ei)_{A,B}$ are affine invariant.
\end{prop}
Furthermore, if the hypersurface lies in a hyperplane which does not contain the origion $O$, since the functions $(ei)_{A,B}$ are affine invariant, we can assume $\mathbf{r}\in\R^n$. Then, we can calculate the invariants by
\begin{gather}
(ei)_{i+1,i+1}=1+\frac{[\mathbf{r}_1-\mathbf{r},\cdots,\mathbf{r}_{ii}-\mathbf{r}_i,\cdots,\mathbf{r}_n-\mathbf{r}]}{[\mathbf{r}_1-\mathbf{r},\cdots,\mathbf{r}_i-\mathbf{r},\cdots,\mathbf{r}_n-\mathbf{r}]},\\
(ei)_{l+1,i+1}=\frac{[\mathbf{r}_1-\mathbf{r},\cdots,\mathbf{r}_{ii}-\mathbf{r}_i,\cdots,\mathbf{r}_n-\mathbf{r}]}{[\mathbf{r}_1-\mathbf{r},\cdots,\mathbf{r}_l-\mathbf{r},\cdots,\mathbf{r}_n-\mathbf{r}]},\\
(ei)_{l+1,j+1}=\frac{[\mathbf{r}_1-\mathbf{r},\cdots,\mathbf{r}_{ij}-\mathbf{r},\cdots,\mathbf{r}_n-\mathbf{r}]}{[\mathbf{r}_1-\mathbf{r},\cdots,\mathbf{r}_l-\mathbf{r},\cdots,\mathbf{r}_n-\mathbf{r}]},
\end{gather}
especially, here $[\dots]$ is the standard determinant in $\R^n$, $i\neq l$ and $i\neq j$.
\section{Sierpinski carpet, Menger sponge and their generalization}
\subsection{Sierpinski carpet}
\begin{figure}[hbtp]
           \centering
            \begin{tabular}{ccc}
           \includegraphics[width=.3\textwidth]{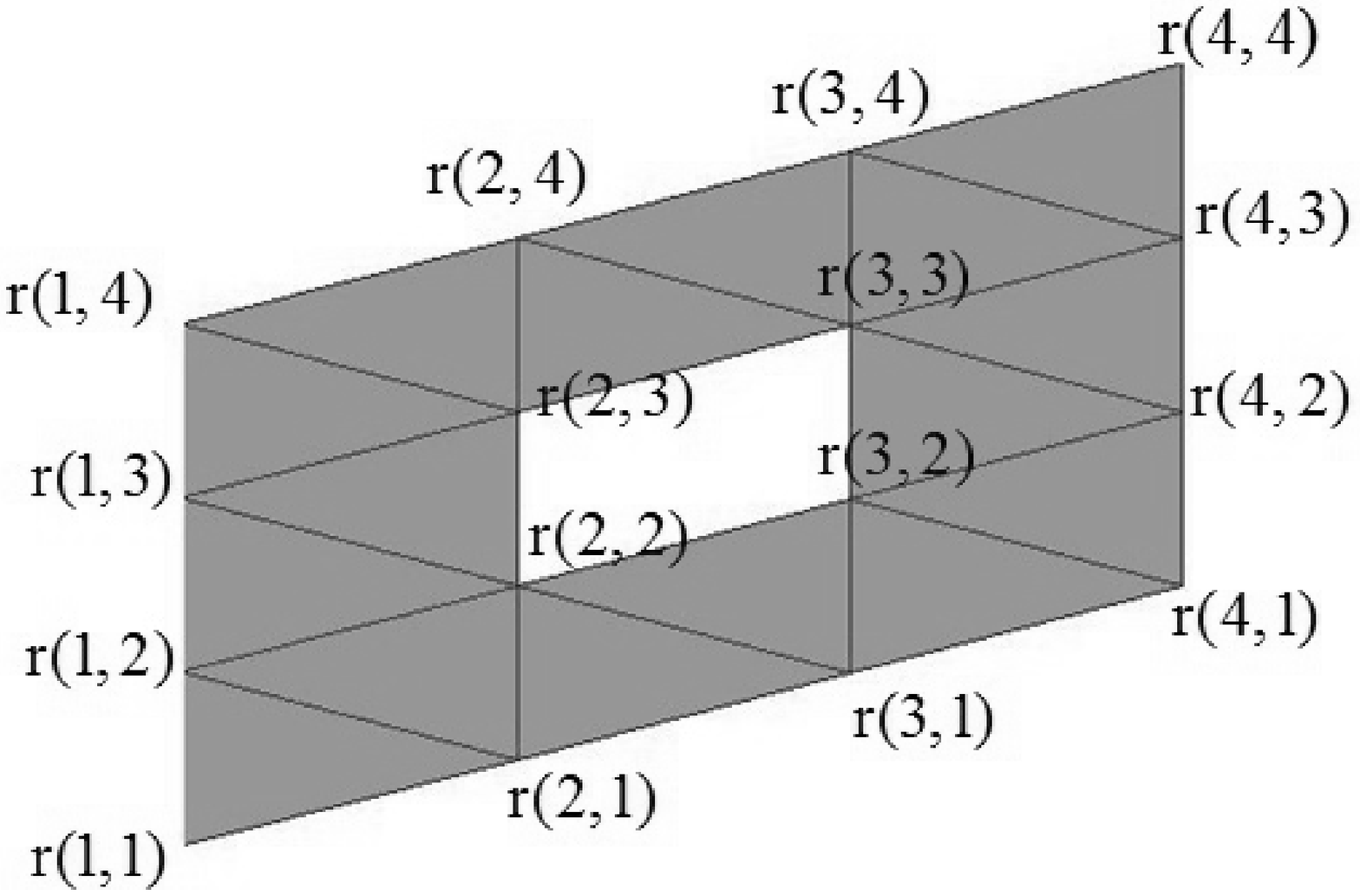}&\includegraphics[width=.3\textwidth]{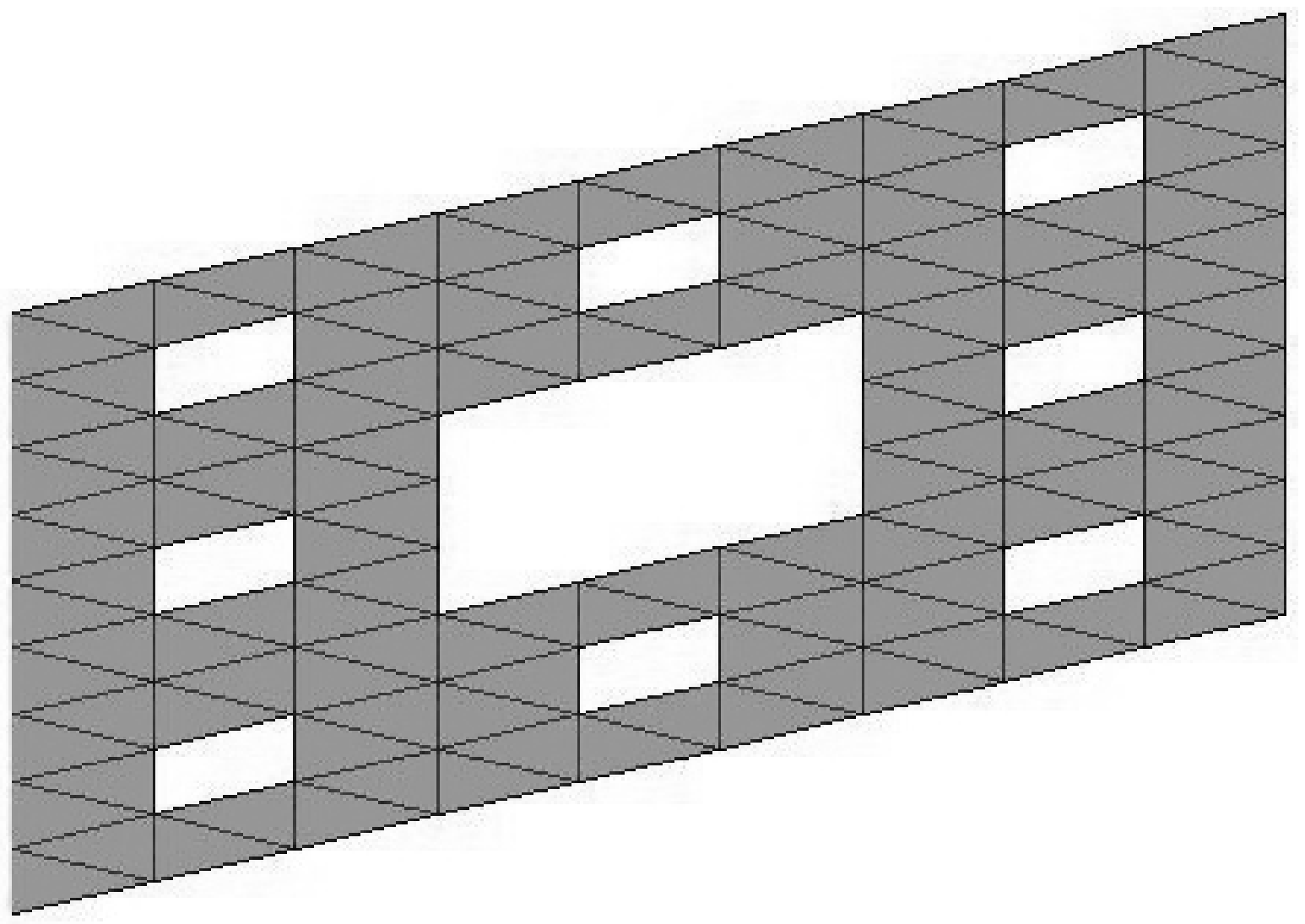}&\includegraphics[width=.3\textwidth]{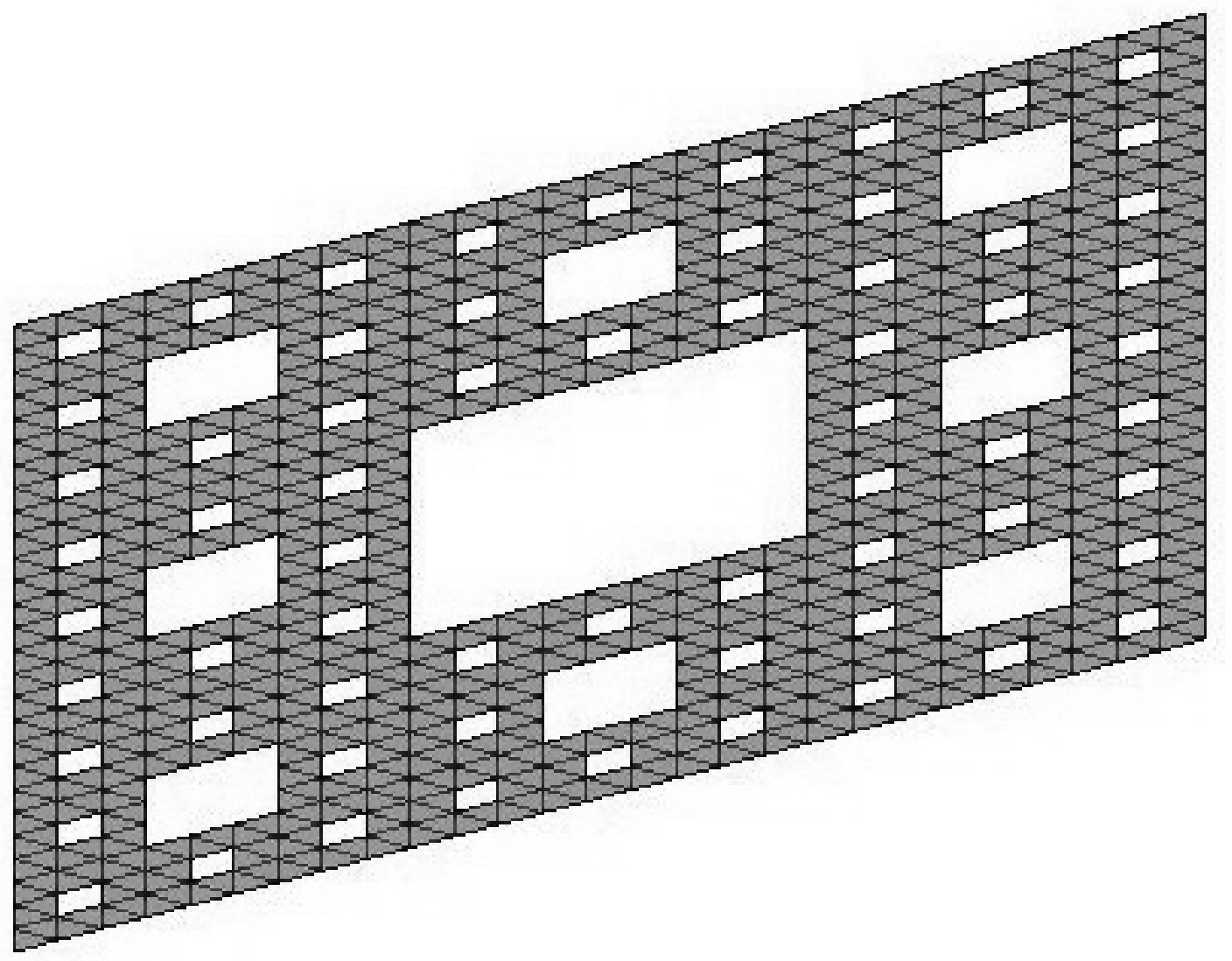}
           \end{tabular}
            \caption{Affine plane Sierpinski carpets and their iterations. }
            \label{fig-carpet-tr}
 \end{figure}

Using the local frame of discrete centro-affine surface in $\R^3$, we can study planar Sierpinski carpet. Since the Sierpinski carpet lie in a plane, all  coefficients in structure equations (\ref{Stru-1-h})-(\ref{Stru-2-h}) are affine invariant. In the left of Figure \ref{fig-carpet-tr}, we mark the point with the vector used in the above section.
\par By a simple calculation, we obtain the structure equation of the vertexes in Figure \ref{fig-carpet-tr}, that is,
\begin{equation}\label{Stru-carpet-St}
  \mathbf{r}_{ij}=\mathbf{r}_i+\mathbf{r}_j-\mathbf{r},\quad i,j=1,2.
\end{equation}
Then it is direct to obtain from Eqs. (\ref{Tran-1-h})-(\ref{Tran-3-h})
\begin{equation}\label{Matrix-carpet-AB}
  M1=\left(
      \begin{array}{ccc}
        0 & -1 & -1 \\
        1 & 2 & 1 \\
        0 & 0 & 1 \\
      \end{array}
    \right),\quad
    M2=\left(
      \begin{array}{ccc}
        0 & -1 & -1 \\
        0 & 1 & 0 \\
        1 & 1 & 2 \\
      \end{array}
    \right),
\end{equation}
which meet the compatibility conditions (\ref{com-con-h}).
\par According to the structure equations, given three non-collinear points $\mathbf{r}(1,1), \mathbf{r}(1,2)$ and
$\mathbf{r}(2,1)$, we can generate all the points in turn. For example, from Eq. (\ref{Stru-carpet-St}), we have
\begin{align}\label{ge-carpet-p}
 \begin{split}
  &\mathbf{r}(2,2)=\mathbf{r}(1,2)+\mathbf{r}(2,1)-\mathbf{r}(1,1),\\
  &\mathbf{r}(1,3)=2\mathbf{r}(1,2)-\mathbf{r}(1,1),\\
  &\mathbf{r}(3,1)=2\mathbf{r}(2,1)-\mathbf{r}(1,1),\\
  &\qquad \cdots \cdots
 \end{split}
\end{align}
The order of these points generated by the structure equation (\ref{Stru-carpet-St}) is
\begin{equation}\label{order-carpet-p}
  \mathbf{r}(1,1), \mathbf{r}(1,2),\mathbf{r}(2,1)\Rightarrow \mathbf{r}(2,2), \mathbf{r}(1,3),\mathbf{r}(3,1),\mathbf{r}(2,3),\mathbf{r}(3,2),\mathbf{r}(3,3),\mathbf{r}(1,4),\mathbf{r}(4,1),\cdots
\end{equation}
Now let us computer which quadrilaterals should be colored. When $n=1$, we obtain a sequence, $$S_1=\{(1,1), (1,2), (2,1),(1,3), (3,1), (2,3), (3,2), (3,3)\},$$
where $(a,b)\in S_1$  corresponds to the a quadrilateral with the four vertexes $\mathbf{r}(a,b), (\mathbf{r}(a+1,b), \mathbf{r}(a+1,b+1),$ and $\mathbf{r}(a,b+1)$. Thus, the sequence $S_1$  implies these eight quadrilaterals should be colored.
\par On the other hand, $\forall (a,b)\in S_1$, $(a,b)$ can be represented by
$$(a,b)=(1,1)+(i,j), \quad i,j=0,1,2,\quad (i,j)\neq (1,1).$$
It is clear that there are $8$ elements in $S_1$.
\par Similarly, when $n=2$,  we also get a sequence $S_2$. As shown in the middle of Figure \ref{fig-carpet-tr}, there are $64$ elements in $S_2$. $\forall (a,b)\in S_2$, $(a,b)$ has the following expression.
$$(a,b)=(1,1)+(i_1,j_1)\times 3^0+(i_2,j_2)\times 3^1, \quad i_k,j_k =0,1,2, \quad (i_k,j_k)\neq (1,1),\quad k=1,2.$$
\par Obviously, at the step $n$, there are $8^n$ colored quadrilaterals. We denote them with a sequence $$S_n=\{(a_1,b_1),\cdots, (a_{8^n},b_{8^n})\}.$$
Exactly, if $(a,b)\in S_n$, $(b,a)\in S_n$. $\forall (a,b)\in S_n$,  $(a,b)$ satisfies that
\begin{equation}\label{Se-carpet}
  (a,b)=(1,1)+(i_1,j_1)\times 3^0
              +(i_2,j_2) \times 3^1+\cdots
              +(i_n,j_n)\times 3^{n-1},
\end{equation}
where
$i_k,j_k=0,1,2, \quad (i_k,j_k)\neq (1,1),\quad k=1,2,\cdots,n$.
\par To sum up, given three non-collinear points and the step $n$, by using Eq. (\ref{ge-carpet-p}) or (\ref{Stru-carpet-St}), after we obtain a point $P$ , it is convenient to verify whether its index $(a,b)$ belongs to $S_n$ in Eq. (\ref{Se-carpet}). If its index $(a,b)$ is in $S_n$, we color the quadrilateral with the vertexes $\mathbf{r}(a,b),\mathbf{r}(a+1,b),\mathbf{r}(a+1,b+1)$ and  $\mathbf{r}(a,b+1)$. Step by step, Sierpinski carpet could be generated. Using the initial points $ \mathbf{r}(1,1)=[0~0], \mathbf{r}(1,2)=[1~1], \mathbf{r}(2,1)=[0~2]$ and $n=1,2,3$, we get three Sierpinski carpets as shown in Figure \ref{fig-carpet-tr}.
Finally, we can answer what is an affine planar Sierpinski carpet.
\begin{prop}\label{prop-carpet-pla}
Exactly the regularities of an affine planar Sierpinski carpet are the structure equations
\begin{equation*}
  \mathbf{r}_{ij}=\mathbf{r}_i+\mathbf{r}_j-\mathbf{r},\quad i,j=1,2,
\end{equation*}
 and the sequence $S_n$ of the colored parallelograms. $\forall (a,b)\in S_n$ corresponds to a quadrilateral with the four vertexes $\mathbf{r}(a,b), (\mathbf{r}(a+1,b), \mathbf{r}(a+1,b+1),$ and $\mathbf{r}(a,b+1)$ and it satisfies
 \begin{equation*}
  (a,b)=(1,1)+(i_1,j_1)\times 3^0
              +(i_2,j_2) \times 3^1+\cdots
              +(i_n,j_n)\times 3^{n-1},
\end{equation*}
where
$i_k,j_k=0,1,2, \quad (i_k,j_k)\neq (1,1),\quad k=1,2,\cdots,n$. There are $8^n$ elements in $S_n$.
\end{prop}
\begin{rem} There is another method to generate a planar Sierpinski carpet. We begin with the indices in the sequence $S_n$ in Eq. (\ref{Se-carpet}), and by using the matrices $M1, M2$ in Eq. (\ref{Matrix-carpet-AB}) and iterative regularities in Eqs. (\ref{Tran-1-h})-(\ref{Tran-3-h}), directly obtain the corresponding points with the indices.  For example, we find the index $(2,3)\in S_n$. Hence, from Eqs. (\ref{Tran-1-h})-(\ref{Tran-3-h}), we obtain
\begin{equation}
  \left(\mathbf{r}(2,3), \mathbf{r}(3,3), \mathbf{r}(2,4)\right)=\left(\mathbf{r}(1,1), \mathbf{r}(2,1), \mathbf{r}(1,2)\right)(M1)(M2)^2.
\end{equation}
These three points $\{\mathbf{r}(2,3), \mathbf{r}(3,3), \mathbf{r}(2,4)\}$ are just we need. Obviously,
$$\left(\mathbf{r}(u,v), \mathbf{r}(u+1,v), \mathbf{r}(u,v+1)\right)=\left(\mathbf{r}(a,b), \mathbf{r}(a+1,b), \mathbf{r}(a,b+1)\right)(M1)^{u-a}(M2)^{v-b}.$$
\end{rem}
\subsection{Menger sponge in $3$ dimensional space}
\begin{figure}[hbtp]
           \centering
            \begin{tabular}{ccc}
           \includegraphics[width=.21\textwidth]{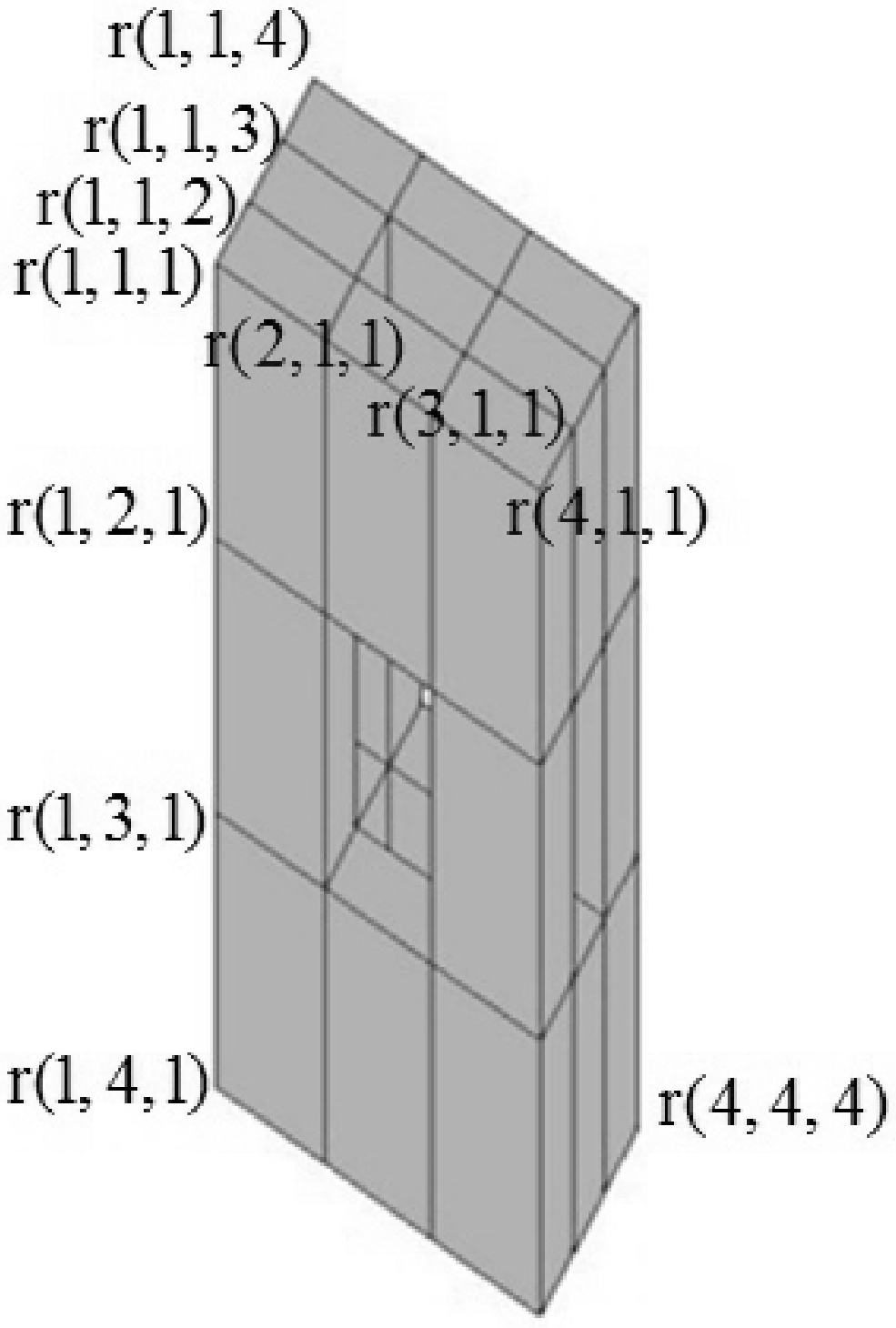}\quad&\quad\includegraphics[width=.12\textwidth]{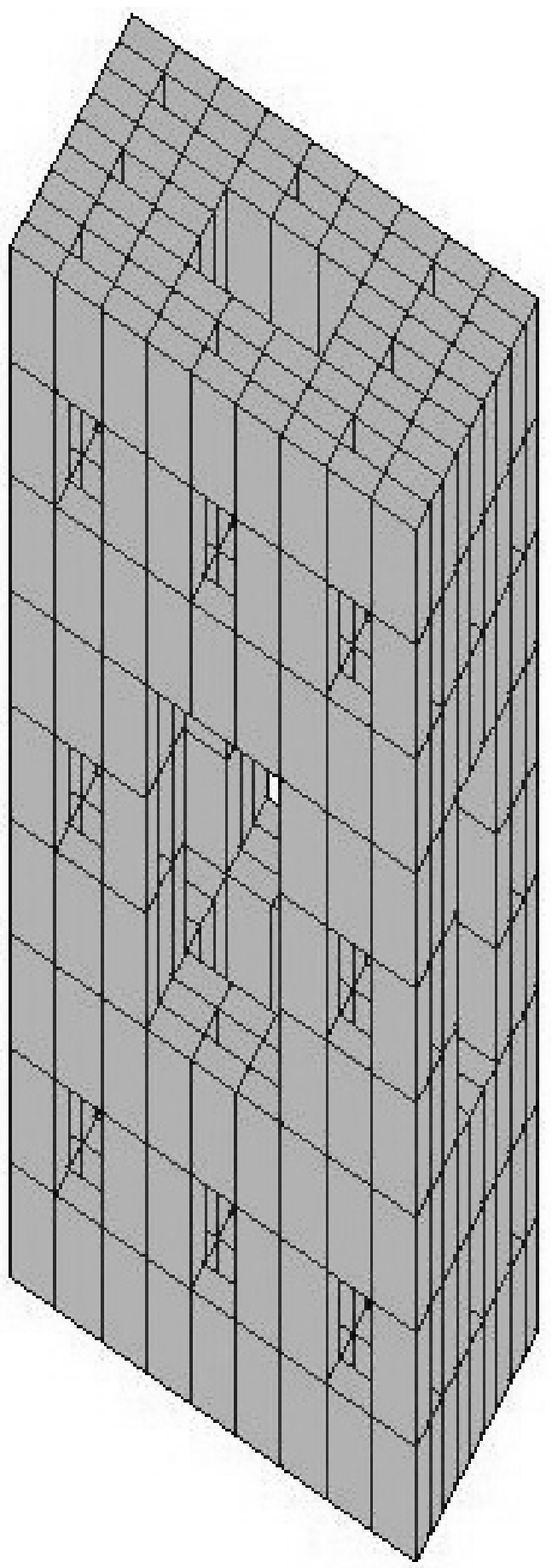}\qquad&\qquad\includegraphics[width=.12\textwidth]{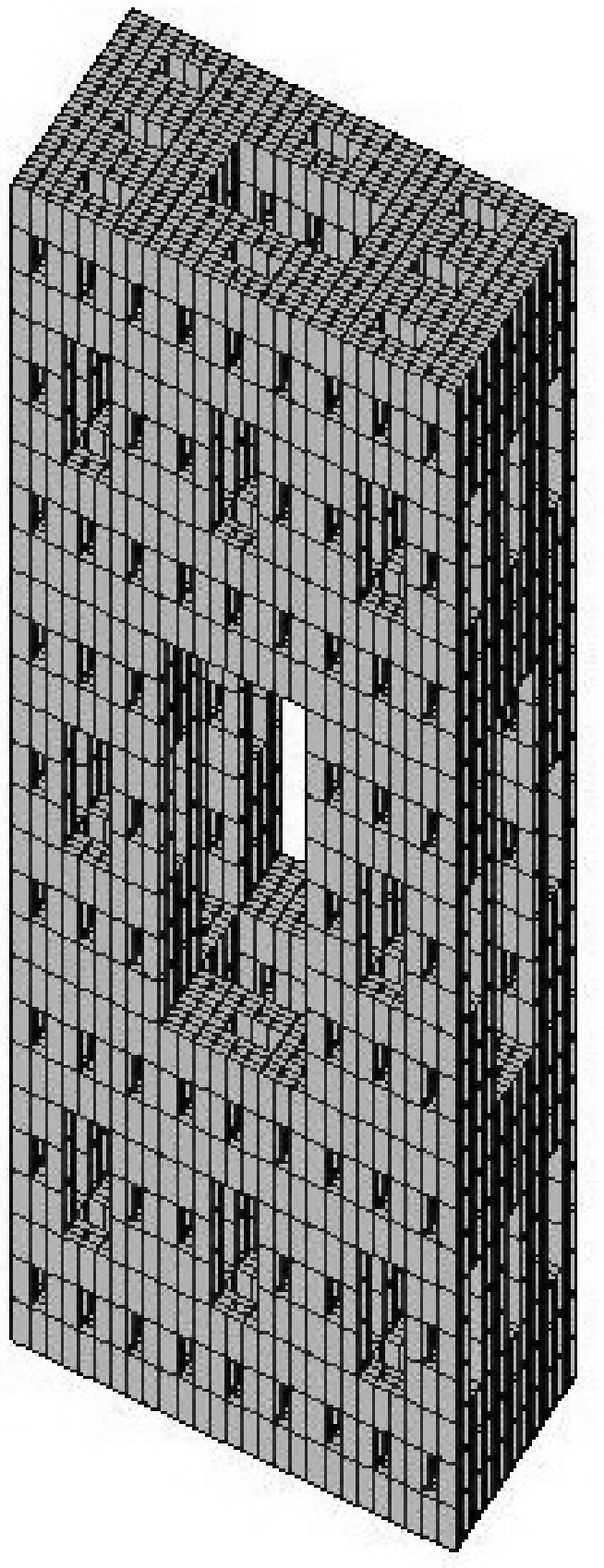}%{ser-sp.eps}
           \end{tabular}
            \caption{Affine spatial Sierpinski carpets with $n=1,2,3$. }
            \label{fig-merge-spm}
 \end{figure}

The three-dimensional version of the Sierpinski carpet is called the Menger sponge. In order to get the affine invariants, we consider it on a hyperplane of $\R^4$. In the left of Figure \ref{fig-merge-spm}, we mark the point with the ordered vector. Hence, it is clear to obtain the relation between these points.
\par By translate it to $\R^3$, using Eqs. (\ref{Stru-1-h})-(\ref{Stru-2-h}), we get the structure equations
\begin{equation}\label{Str-merge-Sp}
  \mathbf{r}_{ij}=\mathbf{r}_i+\mathbf{r}_j-\mathbf{r}, \quad i,j=1,2,3.
\end{equation}
Then according to Eqs. (\ref{Tran-1-h})-(\ref{Tran-3-h}), by a simple calculation, we obtain
\begin{gather}\label{Matrix-merge-Sp}
\begin{split}
M1=\left(
       \begin{array}{cccc}
         0 & -1 & -1 & -1 \\
         1 & 2 & 1 & 1 \\
         0 & 0 & 1 & 0 \\
         0 & 0 & 0 & 1 \\
       \end{array}
     \right)
,\\
M2=\left(
       \begin{array}{cccc}
         0 & -1 & -1 & -1 \\
         0 & 1 & 0 & 0 \\
         1 & 1 & 2 & 1 \\
         0 & 0 & 0 & 1 \\
       \end{array}
     \right)
,\\
M3=\left(
       \begin{array}{cccc}
         0 & -1 & -1 & -1 \\
         0 & 1 & 0 & 0 \\
         0 & 0 & 1 & 0 \\
         1 & 1 & 1 & 2 \\
       \end{array}
     \right).
\end{split}
\end{gather}
It is easy to verify that they satisfy the compatibility conditions (\ref{com-con-h}). Using Corollary \ref{self-sim}, it is obvious that Menger sponge is self-similarity.
\par When $n=1$, we get the sequence $S_1$, and
\begin{align*}S_1=&\{(1,1,1),(2,1,1),(1,2,1),(1,1,2),(3,1,1),(1,3,1),(1,1,3),(3,3,1),(3,1,3),(1,3,3),\\
                  &(2,3,1),(3,2,1),(2,1,3),(3,1,2),(1,2,3),(1,3,2),(2,3,3),(3,2,3),(3,3,2),(3,3,3)\},
\end{align*}
where $\forall (a,b,c)\in S_1$ means a parallelepiped with the vertexes $\mathbf{r}(a,b,c),\mathbf{r}(a+1,b,c), \mathbf{r}(a,b+1,c), \mathbf{r}(a,b,c+1), \mathbf{r}(a+1,b+1,c), \mathbf{r}(a+1,b,c+1), \mathbf{r}(a,b+1,c+1), \mathbf{r}(a+1,b+1,c+1) $  should be colored. There are $20$ elements in $S_1$. Furthermore,
$\forall (a,b,c)\in S_1$, $(a,b,c)$ can be written as
$$(a,b,c)=(1,1,1)+(i,j,k),$$
where $i,j,k=0,1,2$ and any two of $\{i,j,k\}$ can not be $1$  simultaneously.
\par By observing regularities in Figure \ref{fig-merge-spm}, when $n=2$, we obtain the sequence $S_2$, and $\forall (a,b,c)\in S_2$,
$$(a,b,c)=(1,1,1)+(i_1,j_1,k_1)\times 3^0+(i_2,j_2,k_2)\times 3^1,$$
where $i_t,j_t, k_t=0,1,2$ and any two of $\{i_t,j_t,k_t\}$ can not be $1$  simultaneously, $t=1,2$. By a direct calculation,
there are $20^2$ parallelepipeds should be colored.
\par Obviously, at the step $n$, there are $20^n$ colored parallelepipeds. In detail, if $(a,b,c)\in S_n$, $(a,c,b), (b,a,c), (b,c,a), (c,a,b), (c,b,a)\in S_n$, and $(a,b,c)$ could be represented by
\begin{align}\label{Se-merge-p}
  (a,b,c)=(1,1,1)+(i_1,j_1,k_1)\times 3^0+(i_2,j_2,k_2)\times 3^1+\cdots+(i_n,j_n,k_n)\times 3^{n-1},
\end{align}
where $i_t,j_t, k_t=0,1,2$ and any two of $\{i_t,j_t,k_t\}$ can not be $1$  simultaneously, $t=1,2,\cdots,n$.
\par In fact, we also have two methods to generate a Menger sponge by the structure equations.
\par Firstly, we introduce the method using the matrix $M1,M2,M3$ in Eq. (\ref{Matrix-merge-Sp}). Given the initial non-coplanar four points $\mathbf{r}(1,1,1),\mathbf{r}(2,1,1),\mathbf{r}(1,2,1),\mathbf{r}(1,1,2)$. If the index $(a,b,c)\in S_n$, we directly obtain the points by
\begin{align}\label{abc-merge-M}
  &\left(\mathbf{r}(a,b,c),\mathbf{r}(a+1,b,c),\mathbf{r}(a,b+1,c),\mathbf{r}(a,b,c+1)\right)\nonumber\\
  &\qquad\qquad\qquad=\left(\mathbf{r}(1,1,1),\mathbf{r}(2,1,1),\mathbf{r}(1,2,1),\mathbf{r}(1,1,2)\right)(M1)^{a-1}(M2)^{b-1}(M3)^{c-1}.
\end{align}
Therefore, we get the index from Eq. (\ref{Se-merge-p}) in turn, and color the corresponding parallelepiped. In order to reduce computational effort, we also use
\begin{align}\label{abc-merge-M1}
  &\left(\mathbf{r}(u,v,w),\mathbf{r}(u+1,v,w),\mathbf{r}(u,v+1,w),\mathbf{r}(u,v,w+1)\right)\nonumber\\
  &\qquad=\left(\mathbf{r}(a,b,c),\mathbf{r}(a+1,b,c),\mathbf{r}(a,b+1,c),\mathbf{r}(a,b,c+1)\right)(M1)^{u-a}(M2)^{v-b}(M3)^{w-c}.
\end{align}
\par On the other hand, we use Eq. (\ref{Str-merge-Sp}) to generate the points in turn, which implies
\begin{align}
  &\mathbf{r}(u,v,w)=\mathbf{r}(u,v-1,w)+\mathbf{r}(u,v,w-1)-\mathbf{r}(u,v-1,w-1),\\
  &\mathbf{r}(u,v,w)=\mathbf{r}(u-1,v,w)+\mathbf{r}(u,v,w-1)-\mathbf{r}(u-1,v,w-1),\\
  &\mathbf{r}(u,v,w)=\mathbf{r}(u,v-1,w)+\mathbf{r}(u-1,v,w)-\mathbf{r}(u-1,v-1,w),\\
  &\mathbf{r}(u,v,w)=2\mathbf{r}(u-1,v,w)-\mathbf{r}(u-2,v,w),\\
  &\mathbf{r}(u,v,w)=2\mathbf{r}(u,v-1,w)-\mathbf{r}(u,v-2,w),\\
  &\mathbf{r}(u,v,w)=2\mathbf{r}(u,v,w-1)-\mathbf{r}(u,v,w-2).
\end{align}
Using the initial points $ \mathbf{r}(1,1,1)=[0~0~0]$,$\mathbf{r}(2,1,1)=[1~1~0]$,$\mathbf{r}(1,2,1)=[0~2~0]$,$\mathbf{r}(1,1,2)=[0~0~3]$ and $n=1,2,3$, we get three Menger sponges as shown in Figure \ref{fig-serp-spm}.
To sum up, we also have the similar properties for the Menger sponge in $3$ dimensional space as the Sierpinski carpets
\begin{prop}\label{pro-merge-sp}
An affine Menger sponge in $3$-dimensional space can be generated by the structure equations
\begin{equation*}
  \mathbf{r}_{ij}=\mathbf{r}_i+\mathbf{r}_j-\mathbf{r}, \quad i,j=1,2,3,
\end{equation*} and the sequence $S_n$ of parallelepipeds filled with color, $\forall (a,b,c)\in S_n$ satisfies
\begin{align*}
  (a,b,c)=(1,1,1)+(i_1,j_1,k_1)\times 3^0+(i_2,j_2,k_2)\times 3^1+\cdots+(i_n,j_n,k_n)\times 3^{n-1},
\end{align*}
where $i_t,j_t, k_t=0,1,2$ and any two of $\{i_t,j_t,k_t\}$ can not be $1$  simultaneously, $t=1,2,\cdots,n$. $(a,b,c)$ denotes the parallelepiped with the vertexes
$\mathbf{r}(a,b,c), \mathbf{r}(a+1,b,c),\mathbf{r}(a,b+1,c), \mathbf{r}(a,b,c+1), \mathbf{r}(a+1,b+1,c), \mathbf{r}(a+1,b,c+1), \mathbf{r}(a,b+1,c+1),\mathbf{r}(a+1,b+1,c+1).$ There are $20^n$ elements in $S_n$.
\end{prop}
\subsection{Menger sponge in $n$-dimensional space}
Exactly, according to Proposition \ref{prop-carpet-pla} and Proposition \ref{pro-merge-sp}, it is safe to give the definition of an affine Menger sponge in $n$-dimensional space, where $n>1$.
\begin{defn}{\bf(Affine Menger sponge in $n$-dimensional space)} In $n$-dimensional space, where $n>1$,  an affine Menger sponge consists of the structure equations
\begin{equation}\label{n-str}
\mathbf{r}_{ij}=\mathbf{r}_i+\mathbf{r}_j-\mathbf{r}, \quad i,j=1,2,\cdots,n,
\end{equation}
and the sequence $S_m$ of $n$-dimensional parallelehedra filled with color, if $(a_1,a_2,\cdots, a_n)\in S_m$, then $(a_1,a_2,\cdots, a_n)$ could be represented by
\begin{align}
  (a_1,a_2,\cdots, a_n)=(1,1,\cdots,1)+(i_{11},i_{21},\cdots,i_{n1})\times 3^0
              &+(i_{12},i_{22},\cdots,i_{n2})\times 3^1\nonumber\\
              &+\cdots
              +(i_{1m},i_{2m},\cdots,i_{nm})\times 3^{m-1},
\end{align}
where $i_{st}\in \{0,1,2\}$, any two of $\{i_{s1}, i_{s2},\cdots, i_{sn}\}$ can not be $1$ simultaneously. $(a_1,a_2,\cdots, a_n)$ denotes an $n$-dimensional parallelehedra with the vertexes $\mathbf{r}(a_1,a_2,\cdots, a_n)$, $\mathbf{r}(a_1+1,a_2,\cdots, a_n)$, $\mathbf{r}(a_1,a_2+1,\cdots, a_n)$,$\cdots$,$\mathbf{r}(a_1+1,a_2+1,\cdots, a_n+1)$. There are $$\left(3^n
-\left(\begin{array}{c}n \\2 \\\end{array}\right)\times2^{n-2}-\left(\begin{array}{c}n \\3 \\\end{array}\right)\times2^{n-3}-\cdots-\left(\begin{array}{c}n \\ n-1 \\\end{array}\right)\times2-\left(\begin{array}{c}n \\ n \\\end{array}\right)\right)^m$$ in $S_m$, where $\left(\begin{array}{c}n \\r \\\end{array}\right), r=1,2,\cdots, n$, are the binomial coefficients.
\end{defn}
Now, in $4$-dimensional space, we choose the initial points $$\mathbf{r}(1,1,1,1)=(0,0,0,0)^{\mathrm{Tran}},\mathbf{r}(2,1,1,1)=(1,0,0,0)^{\mathrm{Tran}}, \mathbf{r}(1,2,1,1)=(0,1,0,0)^{\mathrm{Tran}},$$ $$\mathbf{r}(1,1,2,1)=(0,0,1,0)^{\mathrm{Tran}},\mathbf{r}(1,1,1,2)=(0,0,0,1)^{\mathrm{Tran}},$$
where the first three are spatial coordinates $(x,y,z)$ and the fourth represents time $t$. It is no difficult to obtain the following matrices from Eqs. (\ref{Tran-1-h})-(\ref{Tran-3-h}), (\ref{n-str})
\begin{gather}
\begin{split}
M1=\left(
       \begin{array}{ccccc}
         0 & -1 & -1 & -1 &-1\\
         1 & 2 & 1 & 1 &1\\
         0 & 0 & 1 & 0 &0\\
         0 & 0 & 0 & 1 &0\\
         0 & 0 & 0 & 0 &1\\
       \end{array}
     \right)
,\quad
M2=\left(
       \begin{array}{ccccc}
         0 & -1 & -1 & -1 &-1\\
         0 & 1 & 0 & 0 &0\\
         1 & 1 & 2 & 1 &1\\
         0 & 0 & 0 & 1 &0\\
         0 & 0 & 0 & 0 &1\\
       \end{array}
     \right)
,\\
M3=\left(
       \begin{array}{ccccc}
         0 & -1 & -1 & -1 &-1\\
         0 & 1 & 0 & 0 &0\\
         0 & 0 & 1 & 0 &0\\
         1 & 1 & 1 & 2 &1\\
         0 & 0 & 0 & 0 &1\\
       \end{array}
     \right),\quad
M4=\left(
       \begin{array}{ccccc}
         0 & -1 & -1 & -1 &-1\\
         0 & 1 & 0 & 0 &0\\
         0 & 0 & 1 & 0 &0\\
         0 & 0 & 0 & 1 &0\\
         1 & 1 & 1 & 1 &2\\
       \end{array}
     \right),
\end{split}
\end{gather}
and
\begin{align}
  &\left(\mathbf{r}(u,v,w,z),\mathbf{r}(u+1,v,w,z),\mathbf{r}(u,v+1,w,z),\mathbf{r}(u,v,w+1,z),\mathbf{r}(u,v,w,z+1)\right)\nonumber\\
  &\qquad=\left(\mathbf{r}(a,b,c,d),\mathbf{r}(a+1,b,c,d),\mathbf{r}(a,b+1,c,d),\mathbf{r}(a,b,c+1,d), \mathbf{r}(a,b,c,d+1)\right)\nonumber\\
  &\qquad \times(M1)^{u-a}(M2)^{v-b}(M3)^{w-c}(M4)^{z-d}.\nonumber
\end{align}
\par It is obvious a $4$-dimensional parallelehedra $(a_1,a_2,a_3,a_4)$ has $16$ vertexes, that is,
$$\mathbf{r}(i_1,i_2,i_3,i_4), i_k\in\{a_k,a_k+1\}, k=1,2,3,4.$$
At the step $n$, there are $48^n$ color $4$-dimensional parallelehedra.
\par
\begin{figure}[hbtp]
           \centering
            \begin{tabular}{c}
           \includegraphics[width=.4\textwidth]{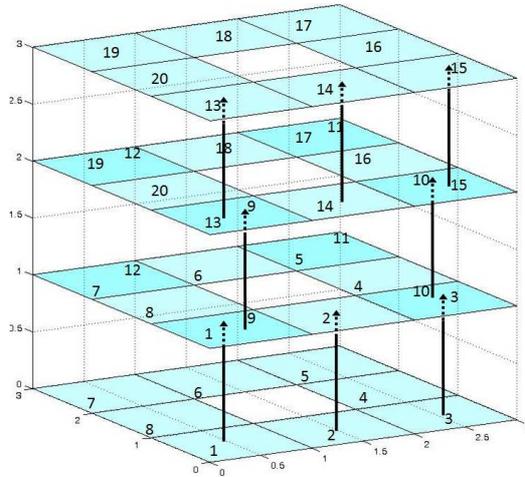}
           \end{tabular}
            \caption{A Menger sponge at $n=1$ from an ant's point of view with $z=0, z=1, z=2,z=3$. }
            \label{fig-carpet-ob1}
 \end{figure}

 In order to understand the Menger sponge in $4$-dimensional space, firstly we observe the Menger sponge in $3$-dimensional space from an ant's point of view.
In this case, $z$-axis can be considered as the axis of time. Hence, at $z=0, z=1,z=2$ and $z=3$, an ant could obtain $4$ Sierpinski carpets as shown in Figure \ref{fig-carpet-ob1}. Exactly, an ant could find the numbers on the surfaces. In fact, in $3$-dimensional space, two squares with the same number are the two bottom surfaces of a cube from $z=t$ to $z=t+1$, $t\in\{0,1,2\}$. At the same place at $z=t$ and $z=t+1$, $t\in\{0,1,2\}$, if the numbers of two squares are not same, then from $z=t$ to $z=t+1$, the place will be empty. If an ant could understand these instructions, it really know what is a Menger sponge in $3$-dimensional space.
 \par
\begin{figure}[hbtp]
           \centering
            \begin{tabular}{cc}
           \includegraphics[width=.45\textwidth]{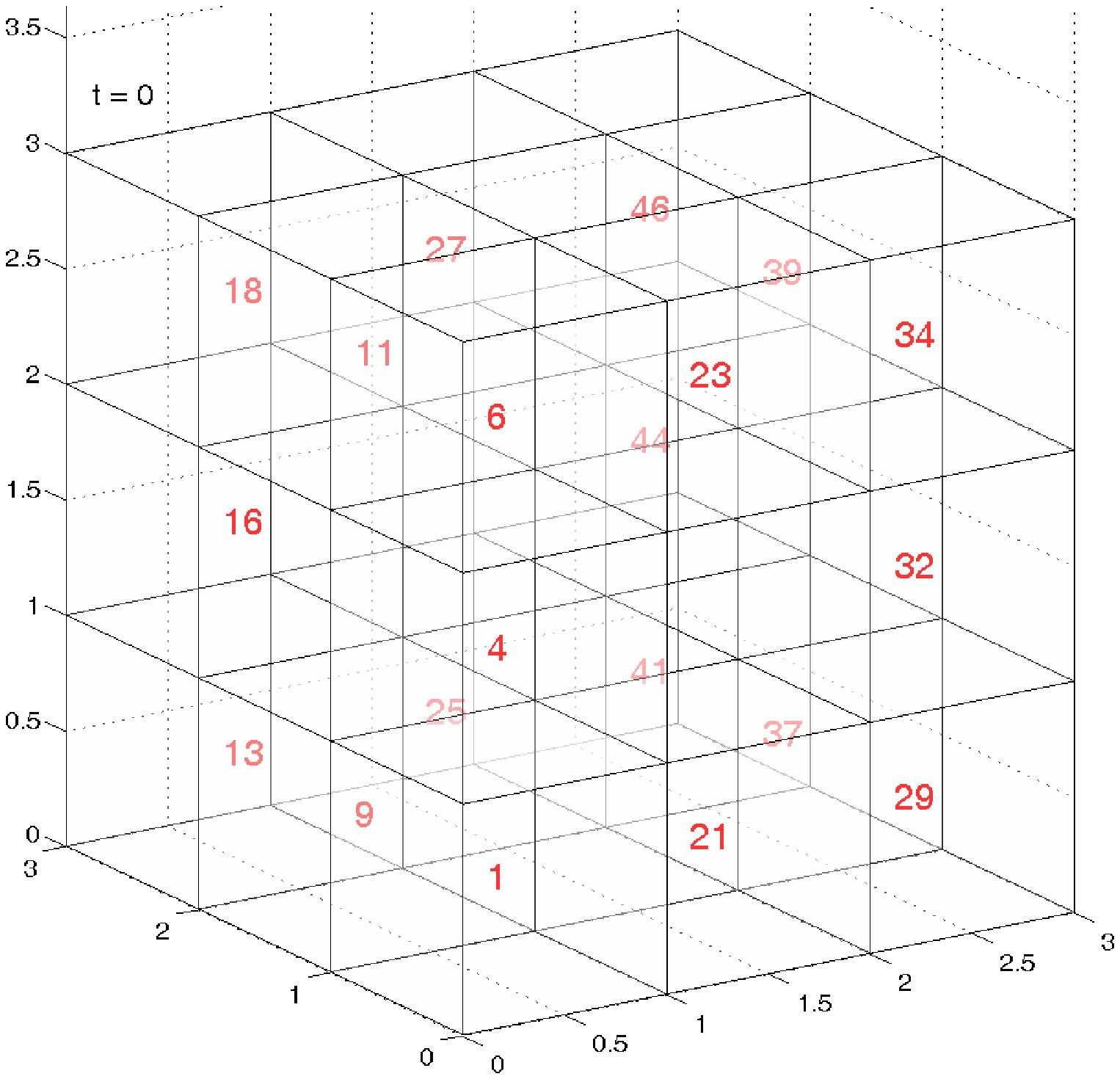}&\includegraphics[width=.45\textwidth]{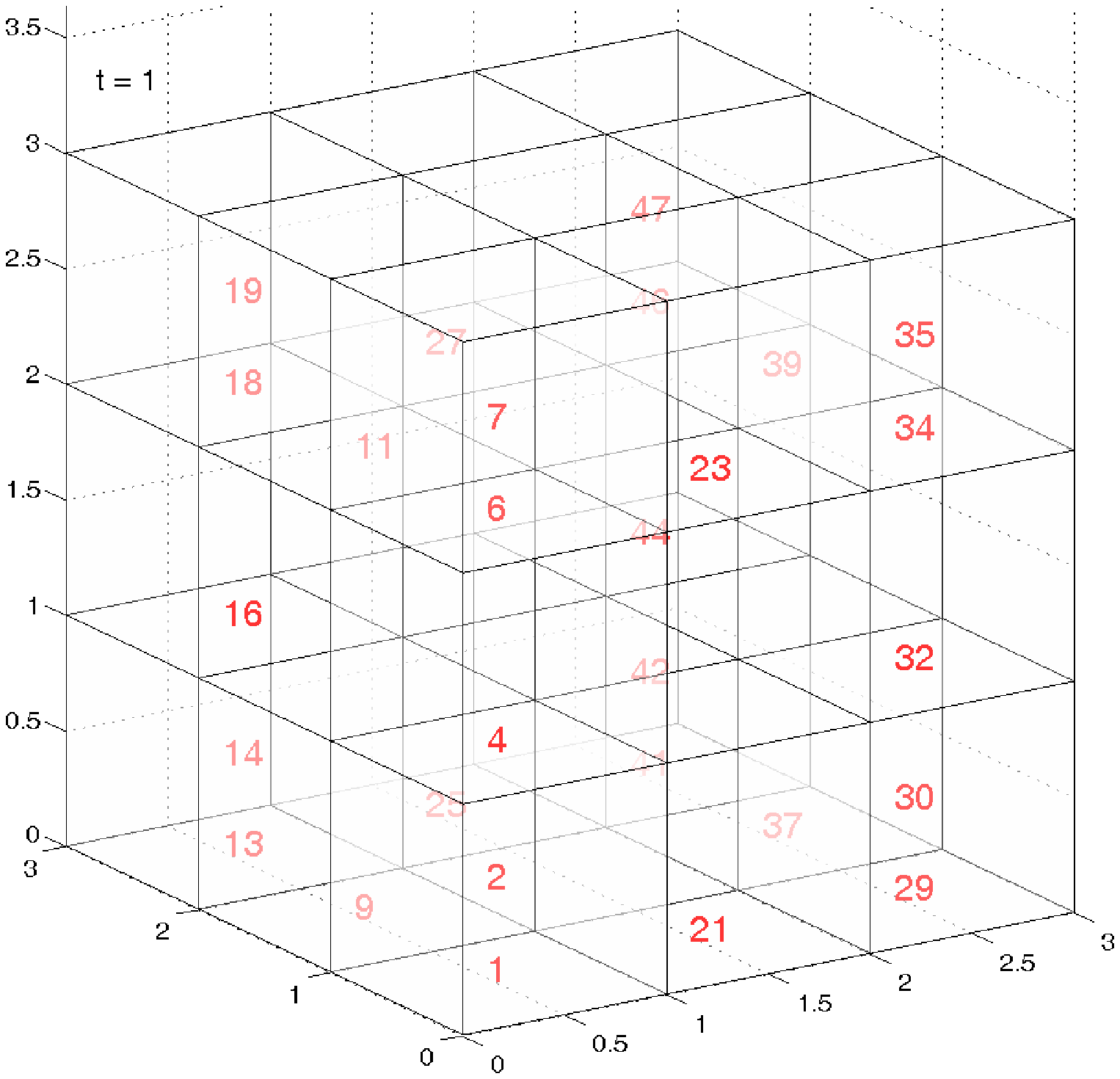}
           \end{tabular}
           \begin{tabular}{ccc}
           \includegraphics[width=.45\textwidth]{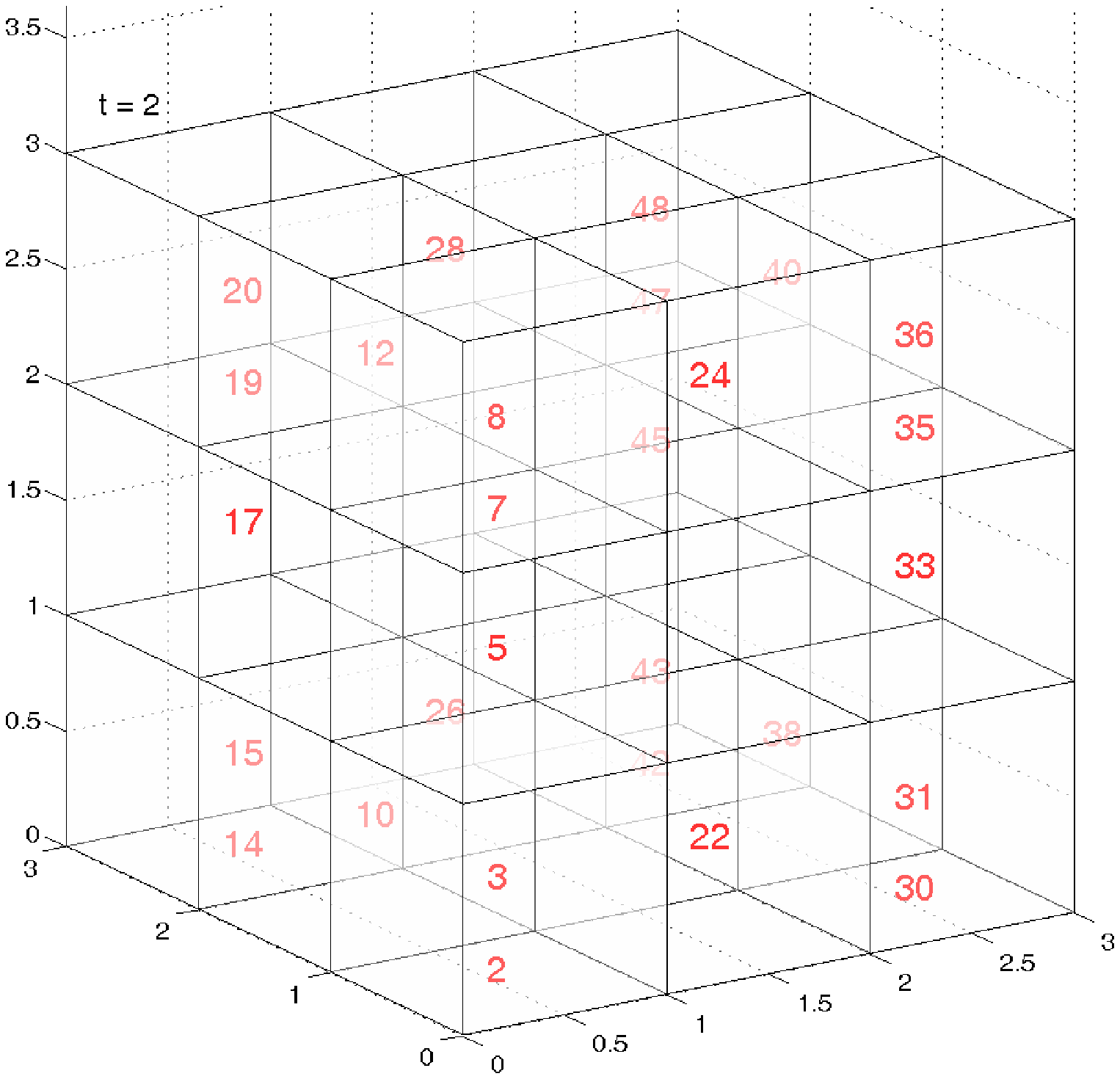}&\includegraphics[width=.45\textwidth]{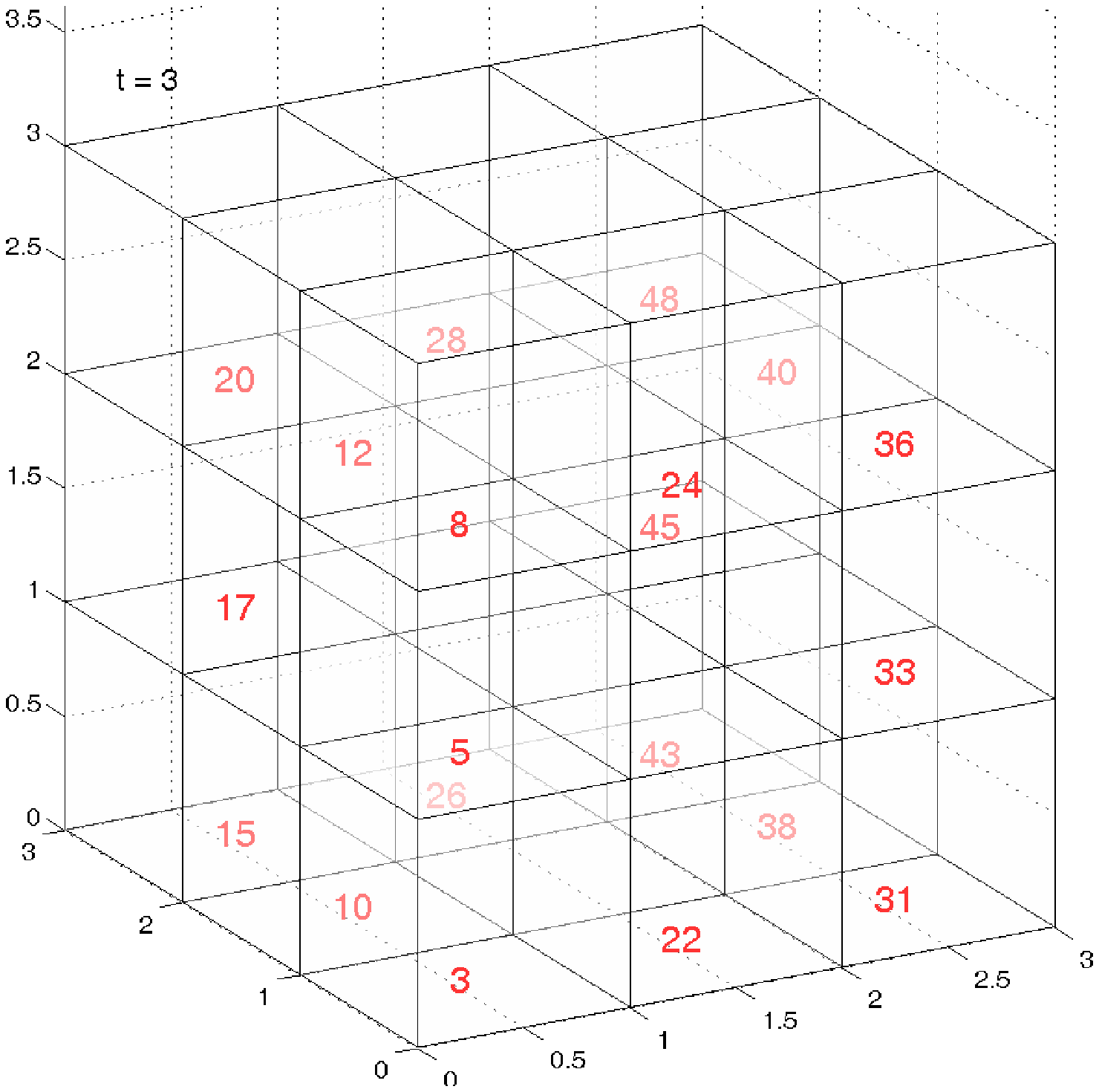}
           \end{tabular}
            \caption{Affine Menger sponge at n=1 in $4$D space with $t=0,1,2,3$. }
            \label{fig-carpet-ob2}
 \end{figure}

 Now let us observe the Monger sponge in $4$-dimensional space.  At $t=0,t=1,t=2$ and $t=3$, we obtain four same Menge sponges in $3$-dimensional space as shown in Figure \ref{fig-carpet-ob2}, but the numbers in them are not all the same. Just as an ant understands a Monger sponge in $3$-dimensional space, we should know, at the same place, if the numbers in two parallelepipeds are same at $t=a$ and $t=a+1$, $a\in\{0,1,2\}$, they are two bottom surfaces of a $4$-dimensional parallelehedra. Else, this place will be empty from $t=a$ to $t=a+1$.
 \section{Sierpinski triangle, Sierpinski pyramid and their generalization}
\subsection{Sierpinski triangle}
\begin{figure}[hbtp]
           \centering
            \begin{tabular}{cc}
           \includegraphics[width=.4\textwidth]{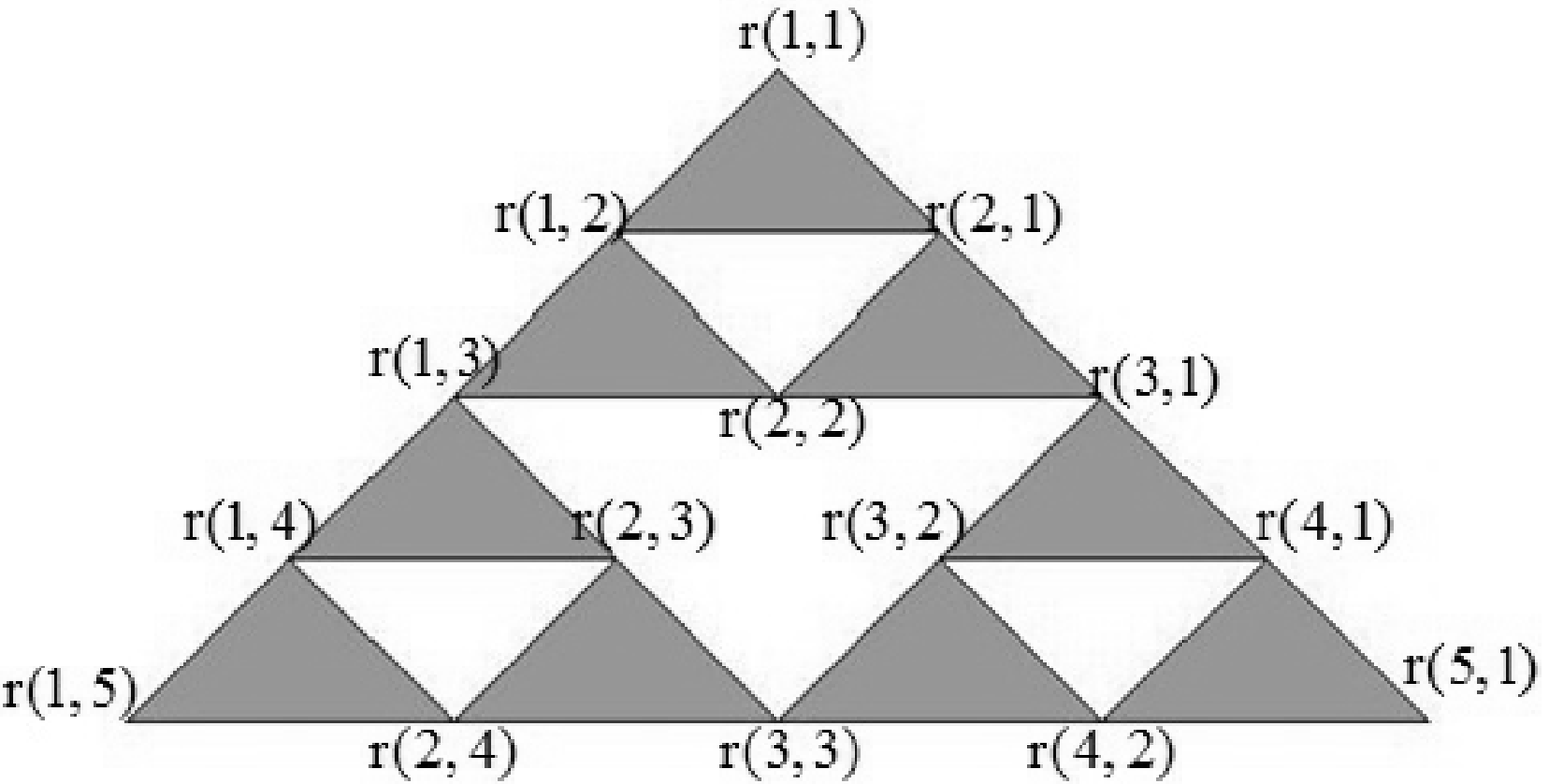}&\includegraphics[width=.3\textwidth]{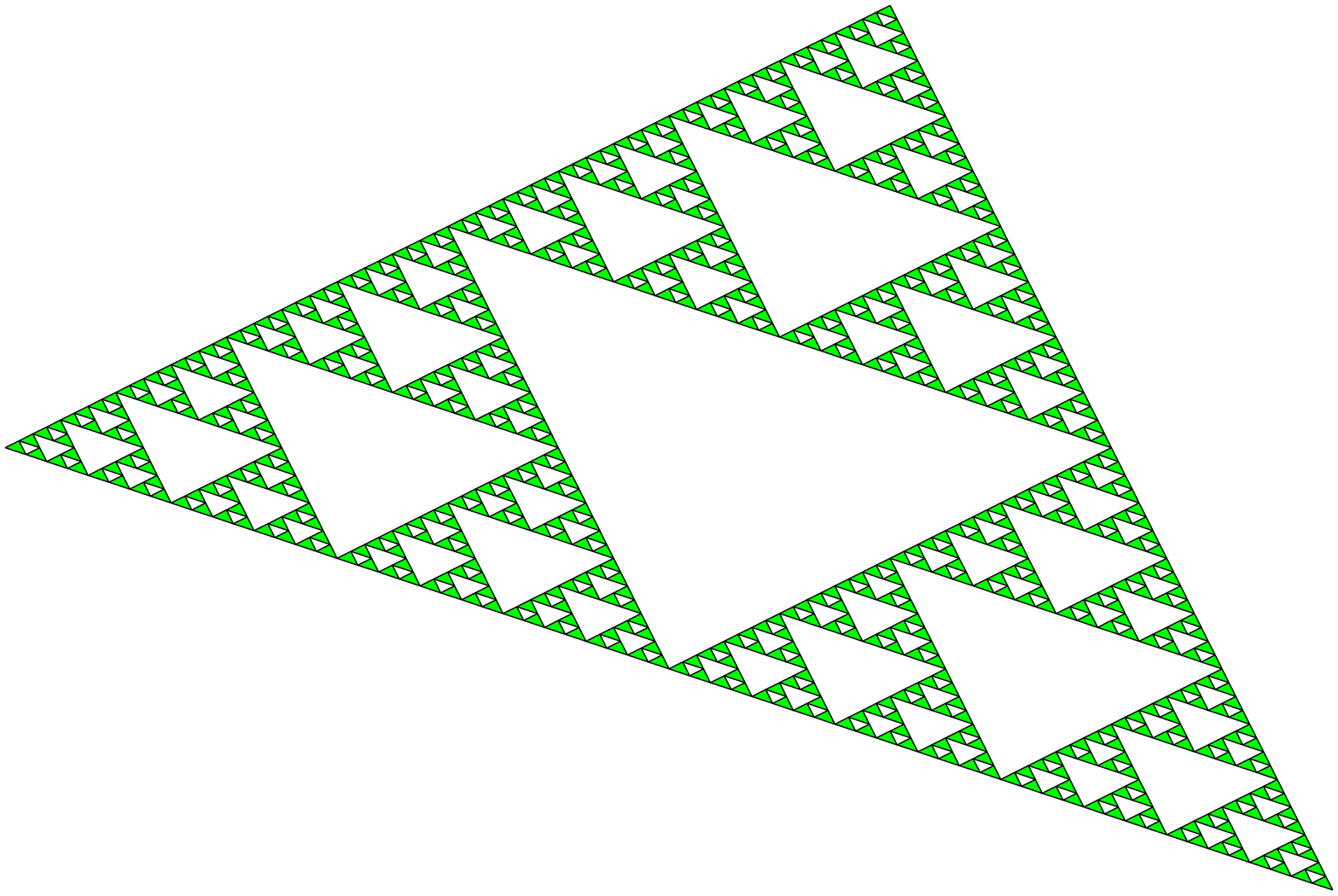}
           \end{tabular}
            \caption{Affine Sierpinski triangle at the step $2$ and $6$. }
            \label{fig-ser-tr}
 \end{figure}

In the left of Figure \ref{fig-ser-tr}, we mark the point with the vector used in the above section.
Similarly, by a simple calculation, we obtain the same structure equations and the matrices as Sierpinski triangle of the vertexes in Figure \ref{fig-ser-tr}, that is,
\begin{equation}\label{Stru-St}
  \mathbf{r}_{ij}=\mathbf{r}_i+\mathbf{r}_j-\mathbf{r},\quad i,j=1,2,
\end{equation}
and
\begin{equation}\label{Matrix-AB}
  M1=\left(
      \begin{array}{ccc}
        0 & -1 & -1 \\
        1 & 2 & 1 \\
        0 & 0 & 1 \\
      \end{array}
    \right),\quad
    M2=\left(
      \begin{array}{ccc}
        0 & -1 & -1 \\
        0 & 1 & 0 \\
        1 & 1 & 2 \\
      \end{array}
    \right).
\end{equation}
Now let us computer which triangles should be colored. When $n=1$, we obtain a sequence, $$(1,1), (1,2), (2,1),$$
where $(1,1)$  corresponds to the triangle $\Delta(\mathbf{r}(1,1),\mathbf{r}(1,2), \mathbf{r}(2,1))$, $(1,2)$ corresponds to the triangle $\Delta(\mathbf{r}(1,2),\mathbf{r}(2,2), \mathbf{r}(1,3))$ and  $(2,1)$ corresponds to the triangle $\Delta(\mathbf{r}(2,1),\mathbf{r}(3,1), \mathbf{r}(2,2))$. It implies these three triangles should be colored.
\par When $n=2$,  we also get a sequence, $$(1,1), (1,2), (2,1),\quad (3,1), (3,2), (4,1),\quad (1,3), (1,4), (2,3).$$
There are $9$ triangles should be colored. The above sequence denote the start point of every triangle, that is, $(m,n)$ represent the triangle with the vertexes $\mathbf{r}(m,n),\mathbf{r}(m,n+1)$, and $\mathbf{r}(m+1,n)$.
Obviously, at the step $n$, there are $3^n$ colored triangles. At the step $n$, if we have the sequence $$S_n=\{(a_1,b_1),\cdots, (a_{3^n},b_{3^n})\},$$
then, at the step $n+1$, the sequence should be
 $$(a_1,b_1),\cdots, (a_{3^n},b_{3^n}),$$
 $$(a_1+2^n,b_1),\cdots, (a_{3^n}+2^n,b_{3^n}),$$
 $$(a_1,b_1+2^n),\cdots, (a_{3^n},b_{3^n}+2^n).$$
Exactly, if $(a,b)\in S_n$, $(b,a)\in S_n$.
\par At the step $n$, if $(a,b)\in S_n$, then $(a,b)$ satisfies that
\begin{equation}\label{Se}
  (a,b)=(1,1)+\left(\begin{array}{c}
                (0,0) \\
                (0,1) \\
                (1,0)
              \end{array}\right)\times 2^0
              +\left(\begin{array}{c}
                (0,0) \\
                (0,1) \\
                (1,0)
              \end{array}\right)\times 2^1+\cdots
              +\left(\begin{array}{c}
                (0,0) \\
                (0,1) \\
                (1,0)
              \end{array}\right)\times 2^{n-1},
\end{equation}
where
$\left(\begin{array}{c}
                (0,0) \\
                (0,1) \\
                (1,0)
              \end{array}\right)$
implies we choose one element from $\{(0,0),(1,0), (0,1)\}$ arbitrarily.
\par On the other hand, we can use matrix to express the iteration of the sequence $S_n$.
$$A_1=\left(
        \begin{array}{cc}
          1 & 1 \\
          1 & 0 \\
        \end{array}
      \right),
      \quad
A_2= \left(
       \begin{array}{cc}
         A_1 & A_1 \\
         A_1 & O \\
       \end{array}
     \right),\cdots
A_{n+1}=\left(
       \begin{array}{cc}
         A_n & A_n \\
         A_n & O \\
       \end{array}
     \right),
$$
where in the matrix $A_n$, the position $(a,b)$ of the element $1$ belongs to $S_n$.
 Using the initial points $ \mathbf{r}(1,1)=[0~0], \mathbf{r}(1,2)=[1~-2], \mathbf{r}(2,1)=[-2~-1]$ and $n=6$, we get Sierpinski triangle as shown in the right of Figure \ref{fig-ser-tr}.
Finally, we can answer what is the Sierpinski triangle.
\begin{prop}
Exactly the regularities of Sierpinski triangle are the structure equation (\ref{Stru-St}) and the sequence $S_n$ in Eq. (\ref{Se}).
\end{prop}

\subsection{Sierpinski pyramid}
\begin{figure}[hbtp]
           \centering
            \begin{tabular}{cc}
           \includegraphics[width=.3\textwidth]{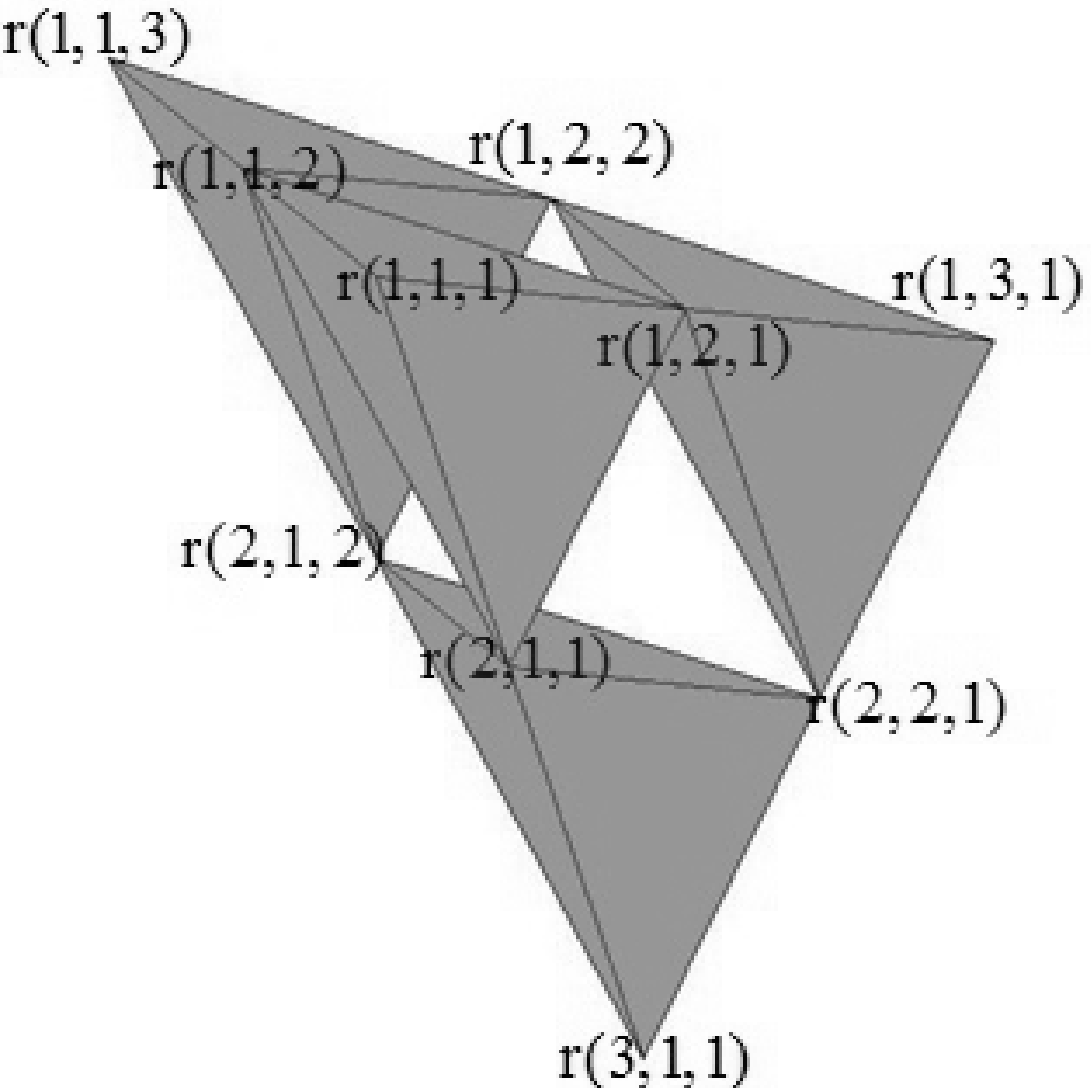}\quad&\quad\includegraphics[width=.3\textwidth]{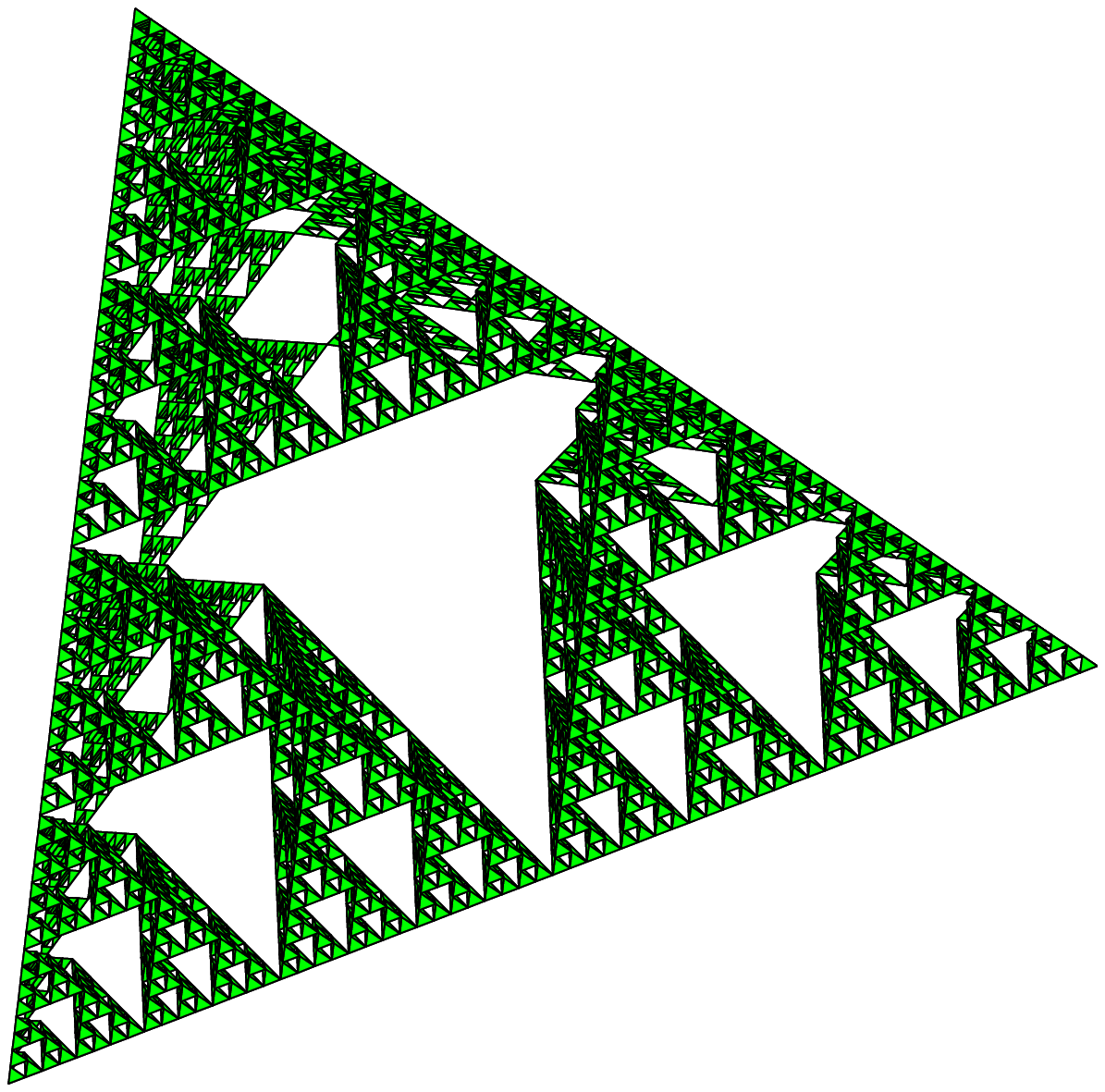}
           \end{tabular}
            \caption{Affine Sierpinski pyramid. }
            \label{fig-serp-spm}
 \end{figure}

In the left of Figure \ref{fig-serp-spm}, we mark the point with the ordered vector.
Similarly, using Eqs. (\ref{Stru-1-h})-(\ref{Stru-2-h}), we get the same structure equations and matrices as Menger sponge.
\begin{equation}\label{Str-Sp}
  \mathbf{r}_{ij}=\mathbf{r}_i+\mathbf{r}_j-\mathbf{r}, \quad i,j=1,2,3,
\end{equation}
and
\begin{gather}\label{Matrix-Sp}
\begin{split}
M1=\left(
       \begin{array}{cccc}
         0 & -1 & -1 & -1 \\
         1 & 2 & 1 & 1 \\
         0 & 0 & 1 & 0 \\
         0 & 0 & 0 & 1 \\
       \end{array}
     \right)
,\\
M2=\left(
       \begin{array}{cccc}
         0 & -1 & -1 & -1 \\
         0 & 1 & 0 & 0 \\
         1 & 1 & 2 & 1 \\
         0 & 0 & 0 & 1 \\
       \end{array}
     \right)
,\\
M3=\left(
       \begin{array}{cccc}
         0 & -1 & -1 & -1 \\
         0 & 1 & 0 & 0 \\
         0 & 0 & 1 & 0 \\
         1 & 1 & 1 & 2 \\
       \end{array}
     \right).
\end{split}
\end{gather}
When $n=1$, we get a sequence, $$(1,1,1), (2,1,1), (1,2,1), (1,1,2)$$ which implies there are  four  tetrahedrons
$$(1,1,1)=(\mathbf{r}(1,1,1),\mathbf{r}(2,1,1), \mathbf{r}(1,2,1), \mathbf{r}(1,1,2)),$$
$$(2,1,1)=(\mathbf{r}(2,1,1),\mathbf{r}(3,1,1), \mathbf{r}(2,2,1),\mathbf{r}(2,1,2)),$$
$$(1,2,1)=(\mathbf{r}(1,2,1),\mathbf{r}(2,2,1), \mathbf{r}(1,3,1),\mathbf{r}(1,2,2)),$$
$$(1,1,2)=(\mathbf{r}(1,1,2),\mathbf{r}(2,1,2), \mathbf{r}(1,2,2),\mathbf{r}(1,1,3)),$$ and they should be colored. Here $(\mathbf{a},\mathbf{b},\mathbf{c},\mathbf{d})$ represents a tetrahedron with the vertexes $\mathbf{a},\mathbf{b},\mathbf{c},\mathbf{d}$.
\par When $n=2$, the sequence can be written as, $$(1,1,1), (2,1,1), (1,2,1), (1,1,2),$$
$$(3,1,1), (4,1,1), (3,2,1), (3,1,2),$$ $$(1,3,1), (2,3,1), (1,4,1), (1,3,2),$$ $$(1,1,3), (2,1,3), (1,2,3), (1,1,4).$$
There are $16$ tetrahedrons should be colored. The above sequences denote the start point of every tetrahedron, that is, $(m,n,l)$ represent a tetrahedron with the vertexes $(\mathbf{r}(m,n,l),\mathbf{r}(m+1,n,l), \mathbf{r}(m,n+1,l), \mathbf{r}(m,n,l+1))$.
\par We have the similar iterative regularities as Sierpinski triangle.
Obviously, at the step $n$, there are $4^n$ colored tetrahedrons.
\par In detail, at the step $n$, if we have the sequence $$S_n=\{(a_1,b_1,c_1),\cdots, (a_{4^n},b_{4^n}, c_{4^n})\},$$
then, at the step $n+1$, the sequence should be
 $$(a_1,b_1,c_1),\cdots, (a_{4^n},b_{4^n,c_{4^n}}),$$
 $$(a_1+2^n,b_1,c_1),\cdots, (a_{4^n}+2^n,b_{4^n},c_{4^n}),$$
 $$(a_1,b_1+2^n,c_1),\cdots, (a_{4^n},b_{4^n}+2^n,c_{4^n}),$$
 $$(a_1,b_1,c_1+2^n),\cdots, (a_{4^n},b_{4^n},c_{4^n}+2^n).$$
Exactly, if $(a,b,c)\in S_n$, $(a,c,b), (b,a,c), (b,c,a), (c,a,b), (c,b,a)\in S_n$.
\par At the step $n$, if $(a,b,c)\in S_n$, then $(a,b,c)$ could be represented by
\begin{align}\label{Se-p}
  (a,b,c)=(1,1,1)+\left(\begin{array}{c}
                (0,0,0) \\
                (1,0,0) \\
                (0,1,0)\\
                (0,0,1)
              \end{array}\right)\times 2^0
              &+\left(\begin{array}{c}
                (0,0,0) \\
                (1,0,0) \\
                (0,1,0)\\
                (0,0,1)
              \end{array}\right)\times 2^1\nonumber\\
              &+\cdots
              +\left(\begin{array}{c}
                (0,0,0) \\
                (1,0,0) \\
                (0,1,0)\\
                (0,0,1)
              \end{array}\right)\times 2^{n-1},
\end{align}
where
$\left(\begin{array}{c}
                (0,0,0) \\
                (1,0,0) \\
                (0,1,0)\\
                (0,0,1)
              \end{array}\right)$
implies we arbitrarily choose one from $\{(0,0,0),(1,0,0), (0,1,0), (0,0,1)\}$.
\par Similarly, we also have two methods to generate Sierpinski pyramid by the structure equations.
Using the initial points $ \mathbf{r}(1,1,1)=[0~0~0]$,$\mathbf{r}(2,1,1)=[-1~-1~-1]$,$\mathbf{r}(1,2,1)=[1~-1~-1]$,$\mathbf{r}(1,1,2)=[1~1~-1]$ and $n=6$, we get Sierpinski triangle as shown in the right of Figure \ref{fig-serp-spm}.
\\ To sum up, we also have the similar properties for Sierpinski pyramid as the Sierpinski triangle
\begin{prop}
The regularities of Sierpinski pyramid are the structure equation (\ref{Str-Sp}) and the sequence $S_n$ in Eq. (\ref{Se-p}).
\end{prop}
In fact, a triangle is a $2$-simplex and a tetrahedron is a $3$-simplex. In the following, we generalize the triangle and the tetrahedron to $n$-simplex.
\subsection{Sierpinski simplex in $n$-dimensional space}
Exactly, according to Eqs. (\ref{Stru-St}), (\ref{Se}), (\ref{Str-Sp}) and (\ref{Se-p}), it is safe to give the definition of affine Sierpinski simplex in $n$-dimensional space,where $n>1$.
\begin{defn}{\bf(Affine Sierpinski simplex)} In $n$-dimensional space,where $n>1$, the structure equations of an affine Sierpinski simplex are
$$\mathbf{r}_{ij}=\mathbf{r}_i+\mathbf{r}_j-\mathbf{r}, \quad i,j=1,2,\cdots,n.$$
At the step $m$, the sequence of Sierpinski simplex is $S_m$, if $(a_1,a_2,\cdots, a_n)\in S_m$, then $(a_1,a_2,\cdots, a_n)$ could be represented by
\begin{align}
  (a_1,a_2,\cdots, a_n)=(1,1,\cdots,1)+&\left(\begin{array}{c}
                (0,0,\cdots,0) \\
                (1,0,\cdots,0) \\
                \cdots\cdots\\
                (0,0,\cdots,1)
              \end{array}\right)\times 2^0
              +\left(\begin{array}{c}
                (0,0,\cdots,0) \\
                (1,0,\cdots,0) \\
                \cdots\cdots\\
                (0,0,\cdots,1)
              \end{array}\right)\times 2^1\nonumber\\
              &+\cdots
              +\left(\begin{array}{c}
                (0,0,\cdots,0) \\
                (1,0,\cdots,0) \\
                \cdots\cdots\\
                (0,0,\cdots,1)
              \end{array}\right)\times 2^{m-1},
\end{align}
where
$\left(\begin{array}{c}
                (0,0,\cdots,0) \\
                (1,0,\cdots,0) \\
                \cdots\cdots\\
                (0,0,\cdots,1)
              \end{array}\right)$
implies to choose arbitrary one from $$\{(0,0,\cdots,0),(1,0,\cdots,0), \cdots, (0,0,\cdots,1)\},$$
and $(a_1,a_2,\cdots, a_n)$ denotes we should color the simplex with the vertexes $\mathbf{r}(a_1,a_2,\cdots, a_n)$, $\mathbf{r}(a_1+1,a_2,\cdots, a_n)$, $\mathbf{r}(a_1,a_2+1,\cdots, a_n)$,$\cdots$,$\mathbf{r}(a_1,a_2,\cdots, a_n+1)$.
\end{defn}
Now, in $4$-dimensional space, we choose the initial points $$\mathbf{r}(1,1,1,1)=(0,0,0,0)^{\mathrm{Tran}},\mathbf{r}(2,1,1,1)=(1,0,0,0)^{\mathrm{Tran}}, \mathbf{r}(1,2,1,1)=(0,1,0,0)^{\mathrm{Tran}},$$ $$\mathbf{r}(1,1,2,1)=(0,0,1,0)^{\mathrm{Tran}},\mathbf{r}(1,1,1,2)=(0,0,0,1)^{\mathrm{Tran}},$$
where the first three are spatial coordinates $(x,y,z)$ and the fourth represents time $t$. It is no difficult to obtain the matrices
\begin{gather}
\begin{split}
M1=\left(
       \begin{array}{ccccc}
         0 & -1 & -1 & -1 &-1\\
         1 & 2 & 1 & 1 &1\\
         0 & 0 & 1 & 0 &0\\
         0 & 0 & 0 & 1 &0\\
         0 & 0 & 0 & 0 &1\\
       \end{array}
     \right)
,\quad
M2=\left(
       \begin{array}{ccccc}
         0 & -1 & -1 & -1 &-1\\
         0 & 1 & 0 & 0 &0\\
         1 & 1 & 2 & 1 &1\\
         0 & 0 & 0 & 1 &0\\
         0 & 0 & 0 & 0 &1\\
       \end{array}
     \right)
,\\
M3=\left(
       \begin{array}{ccccc}
         0 & -1 & -1 & -1 &-1\\
         0 & 1 & 0 & 0 &0\\
         0 & 0 & 1 & 0 &0\\
         1 & 1 & 1 & 2 &1\\
         0 & 0 & 0 & 0 &1\\
       \end{array}
     \right),\quad
M4=\left(
       \begin{array}{ccccc}
         0 & -1 & -1 & -1 &-1\\
         0 & 1 & 0 & 0 &0\\
         0 & 0 & 1 & 0 &0\\
         0 & 0 & 0 & 1 &0\\
         1 & 1 & 1 & 1 &2\\
       \end{array}
     \right),
\end{split}
\end{gather}
and
\begin{align}
  &\left(\mathbf{r}(u,v,w,z),\mathbf{r}(u+1,v,w,z),\mathbf{r}(u,v+1,w,z),\mathbf{r}(u,v,w+1,z),\mathbf{r}(u,v,w,z+1)\right)\nonumber\\
  &\qquad=\left(\mathbf{r}(a,b,c,d),\mathbf{r}(a+1,b,c,d),\mathbf{r}(a,b+1,c,d),\mathbf{r}(a,b,c+1,d), \mathbf{r}(a,b,c,d+1)\right)\nonumber\\
  &\qquad \times(M1)^{u-a}(M2)^{v-b}(M3)^{w-c}(M4)^{z-d}.\nonumber
\end{align}
\par
\begin{figure}[hbtp]
           \centering
            \begin{tabular}{c}
           \includegraphics[width=.8\textwidth]{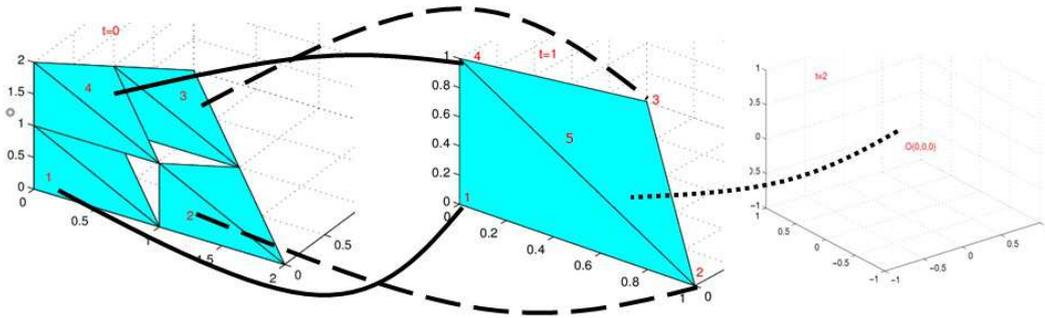}
           \end{tabular}
            \caption{Affine Sierpinski simplex at n=1 in $4$D space with $t=0, t=1, t=2$. }
            \label{fig-serp-ob1}
 \end{figure}

 Then, $n=1$, we obtain the graph as shown in Figure \ref{fig-serp-ob1}.  There are three graphs, which are divided according to $t=0,t=1$ and $t=2$. Of course, they all link together in $4$ dimensional space. In Figure \ref{fig-serp-ob1}, we connect them with the curves, which implies they are in a $5$-simplex with that five vertexes. There is not difference between the full curves and broken curves, only for clarity of the figure. On the other hand, the same numbers also indicate they are in $5$-simplex. For two same numbers, the number in front represents a  tetrahedron, and the last one denotes a point.
\par
\begin{figure}[hbtp]
           \centering
            \begin{tabular}{cc}
           \includegraphics[width=.4\textwidth]{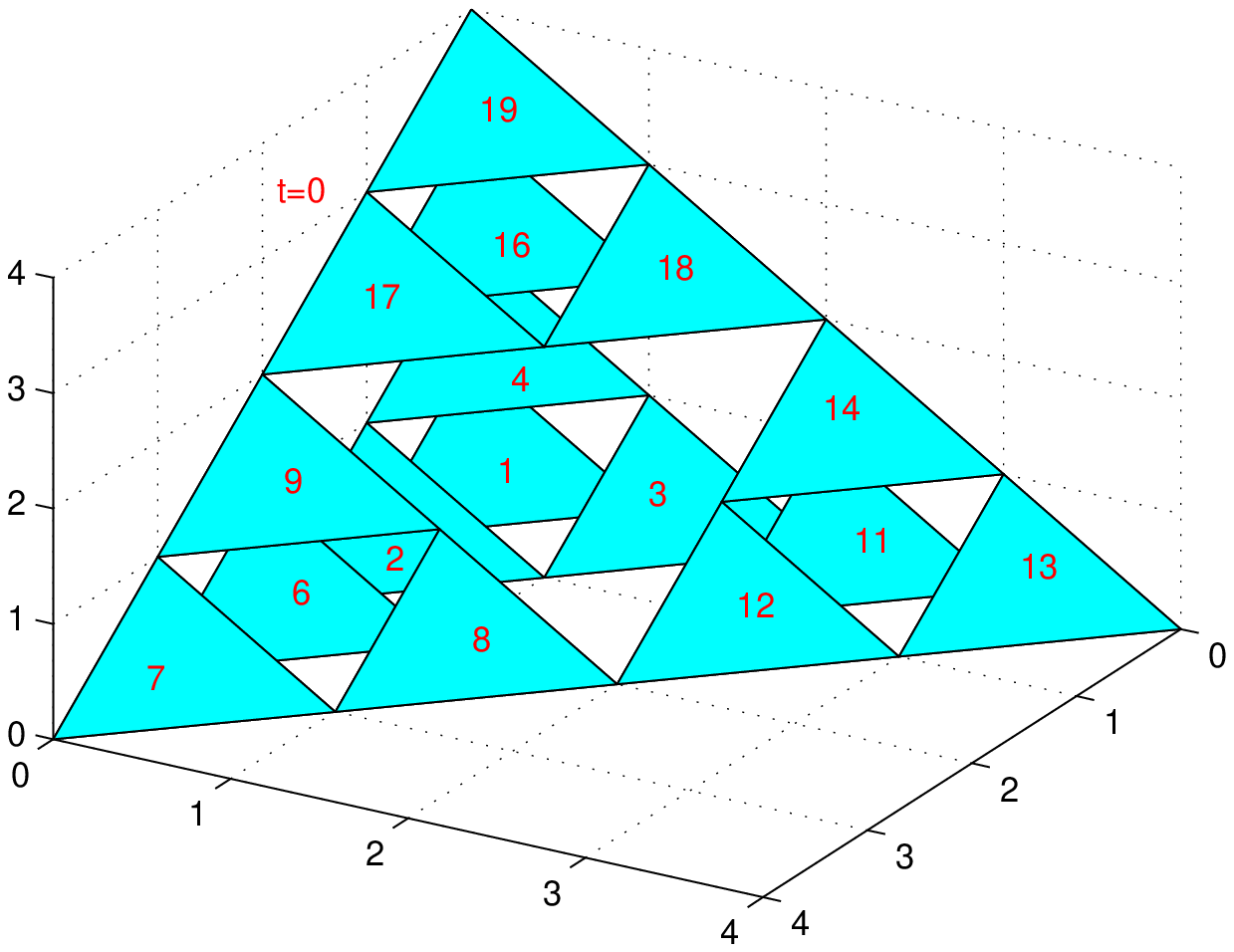}&\includegraphics[width=.4\textwidth]{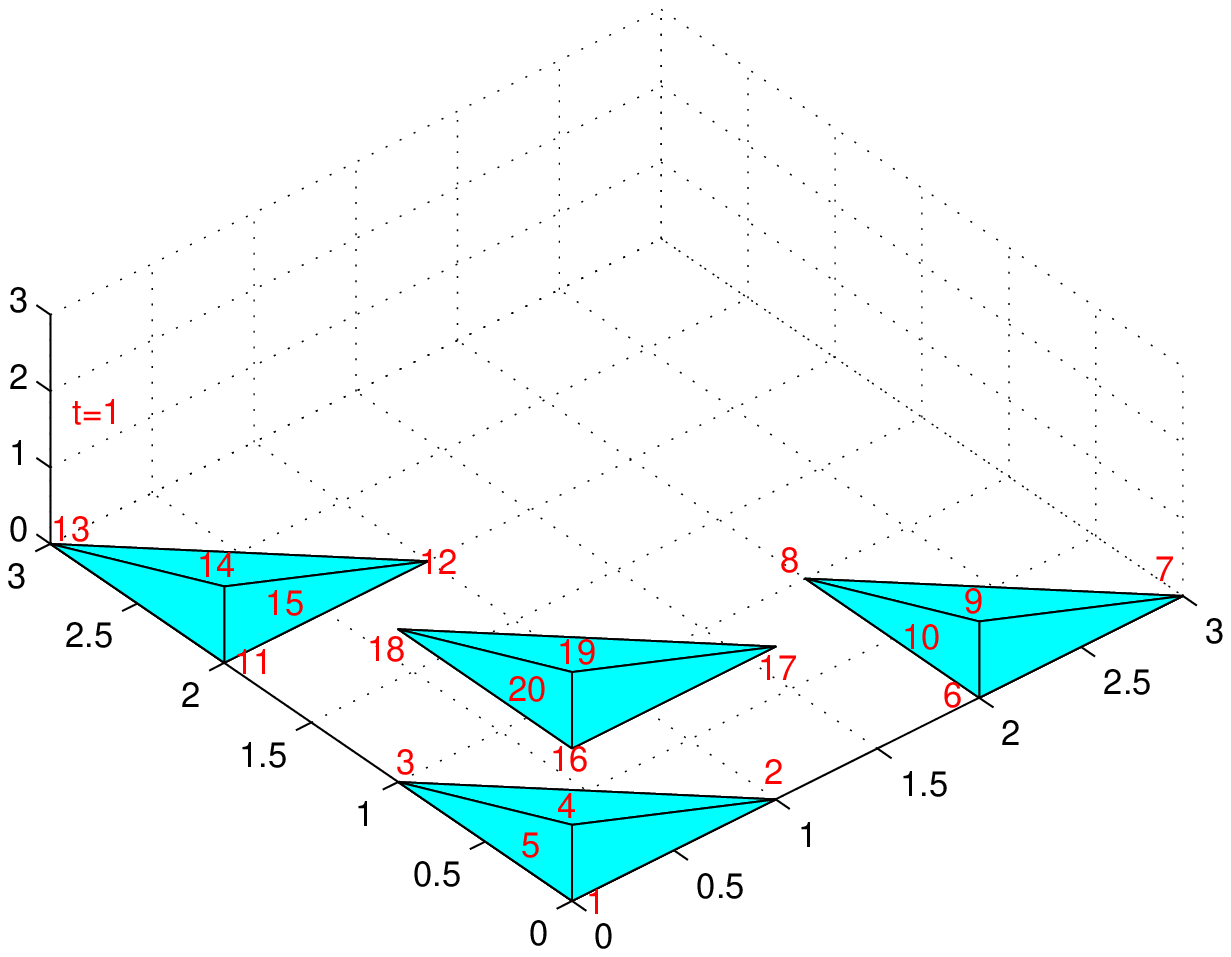}
           \end{tabular}
           \begin{tabular}{ccc}
           \includegraphics[width=.3\textwidth]{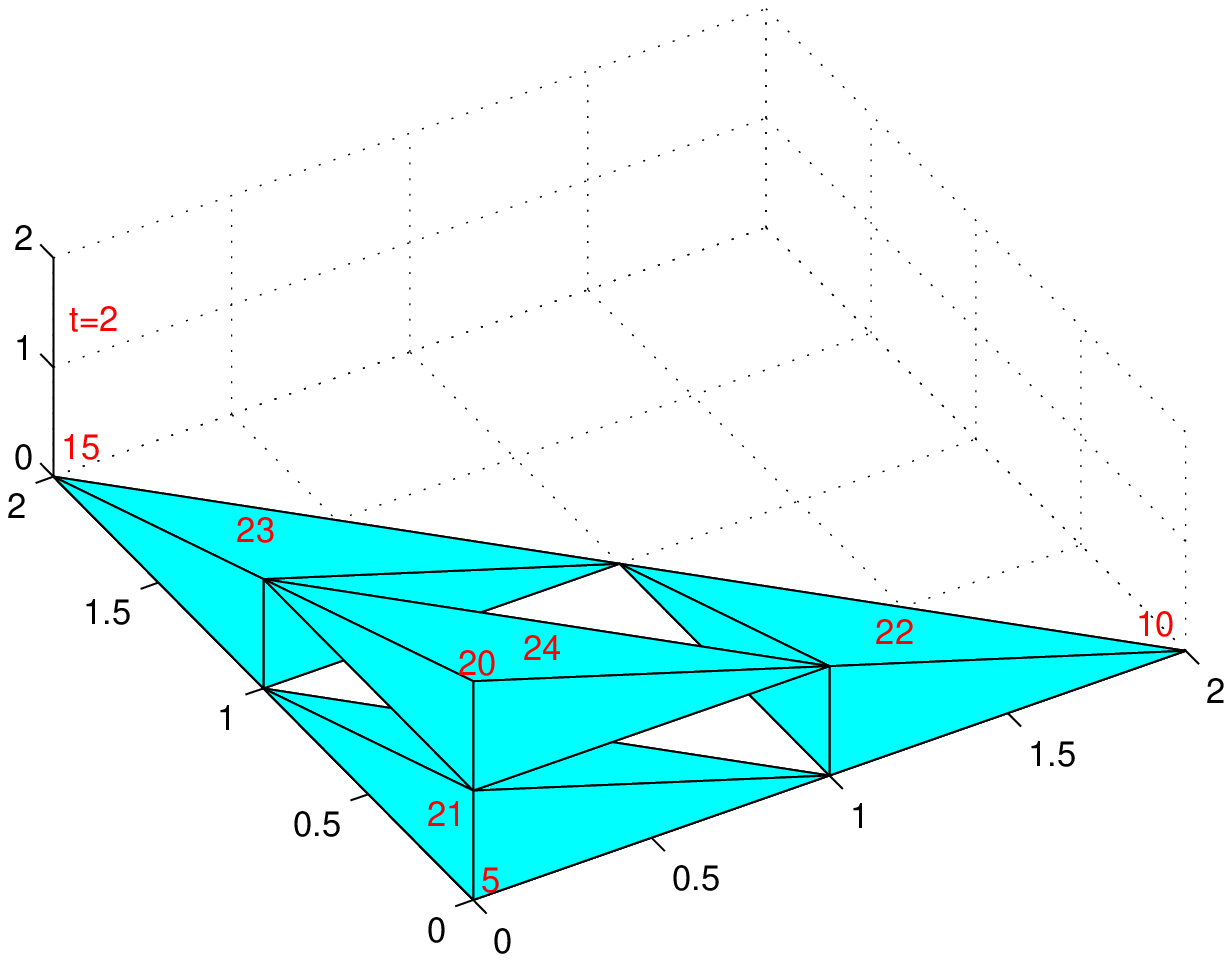}&\includegraphics[width=.3\textwidth]{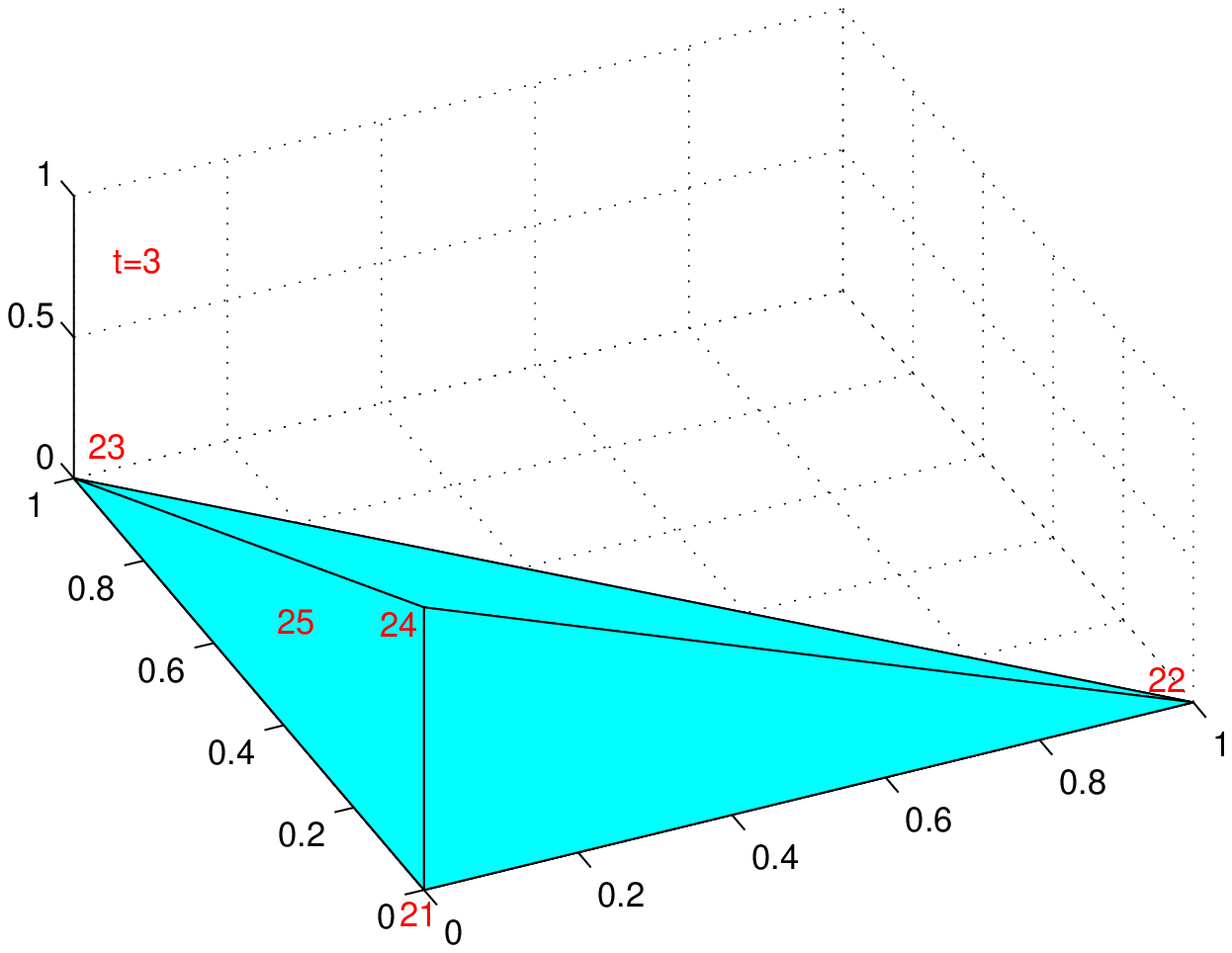}&\includegraphics[width=.3\textwidth]{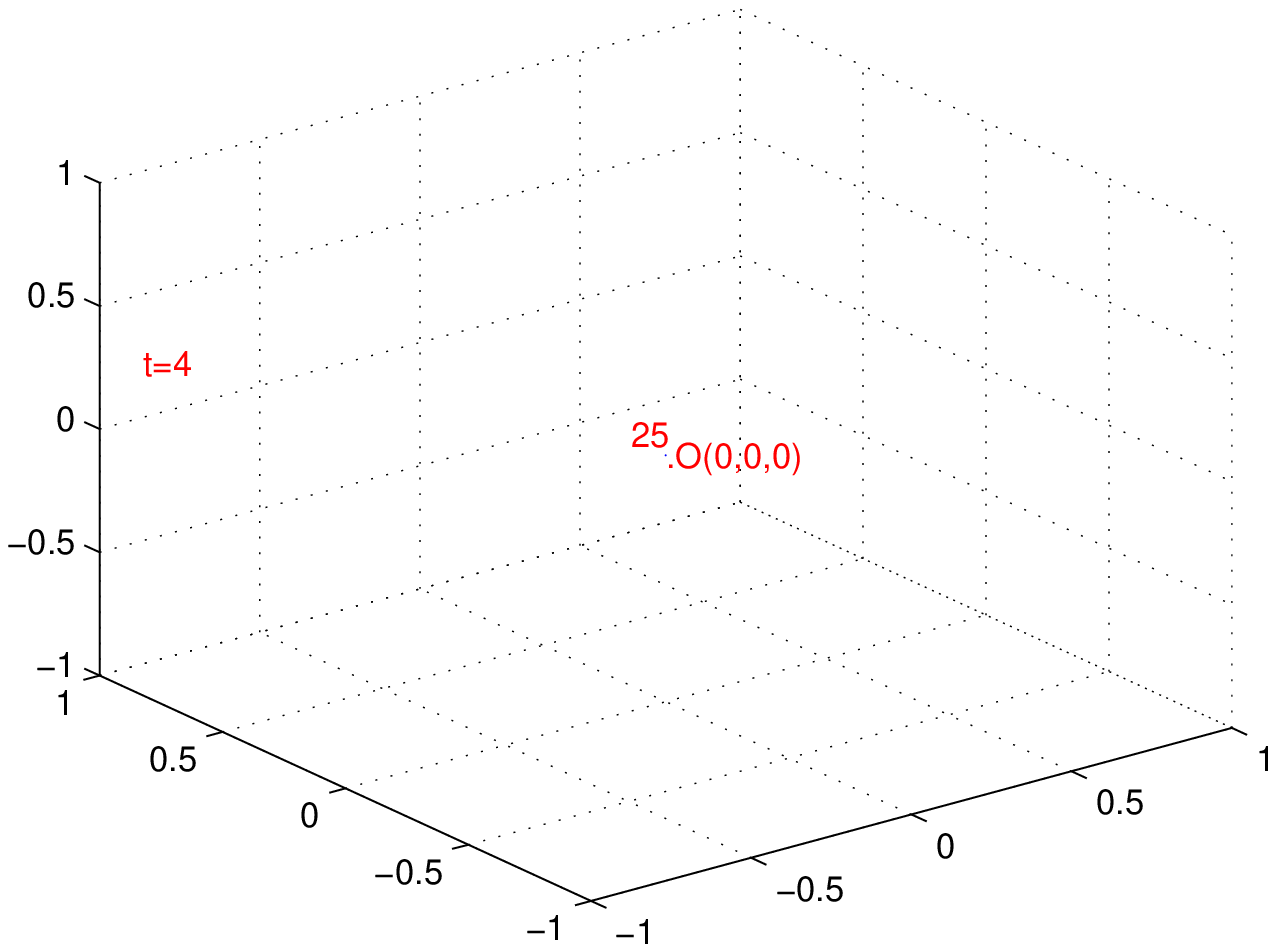}
           \end{tabular}
            \caption{Affine Sierpinski simplex at n=2 in $4$D space with $t=0,1,2,3,4$. }
            \label{fig-serp-ob2}
 \end{figure}

 When $n=2$, it is shown in Figure \ref{fig-serp-ob2}. In this figure, we do not use the curve to connect them. Exactly, we can understand it by the numbers. Of course, the numbers have the same meaning as in Figure \ref{fig-serp-ob1}.
 \par
\begin{figure}[hbtp]
           \centering
            \begin{tabular}{ccc}
           \includegraphics[width=.3\textwidth]{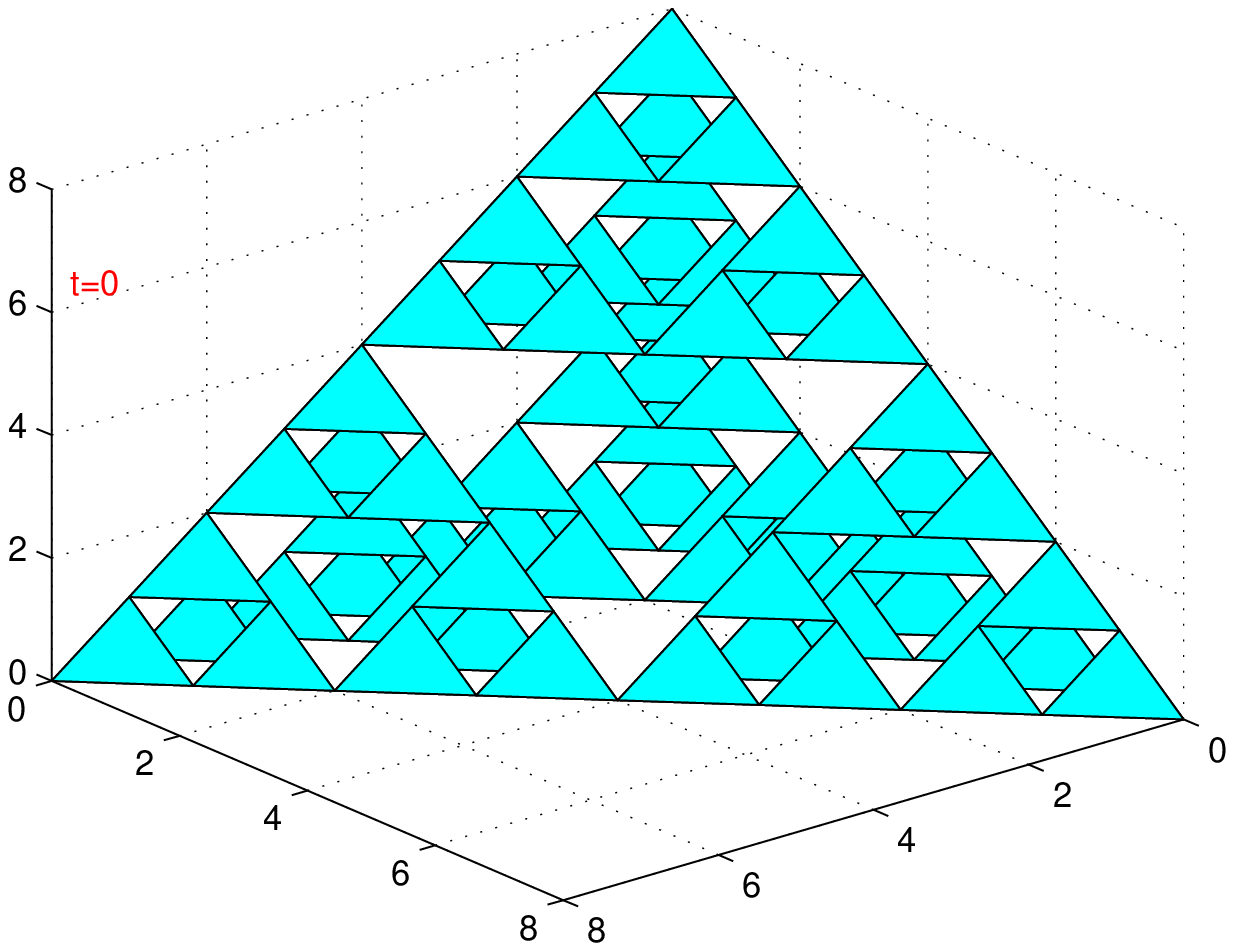}&\includegraphics[width=.3\textwidth]{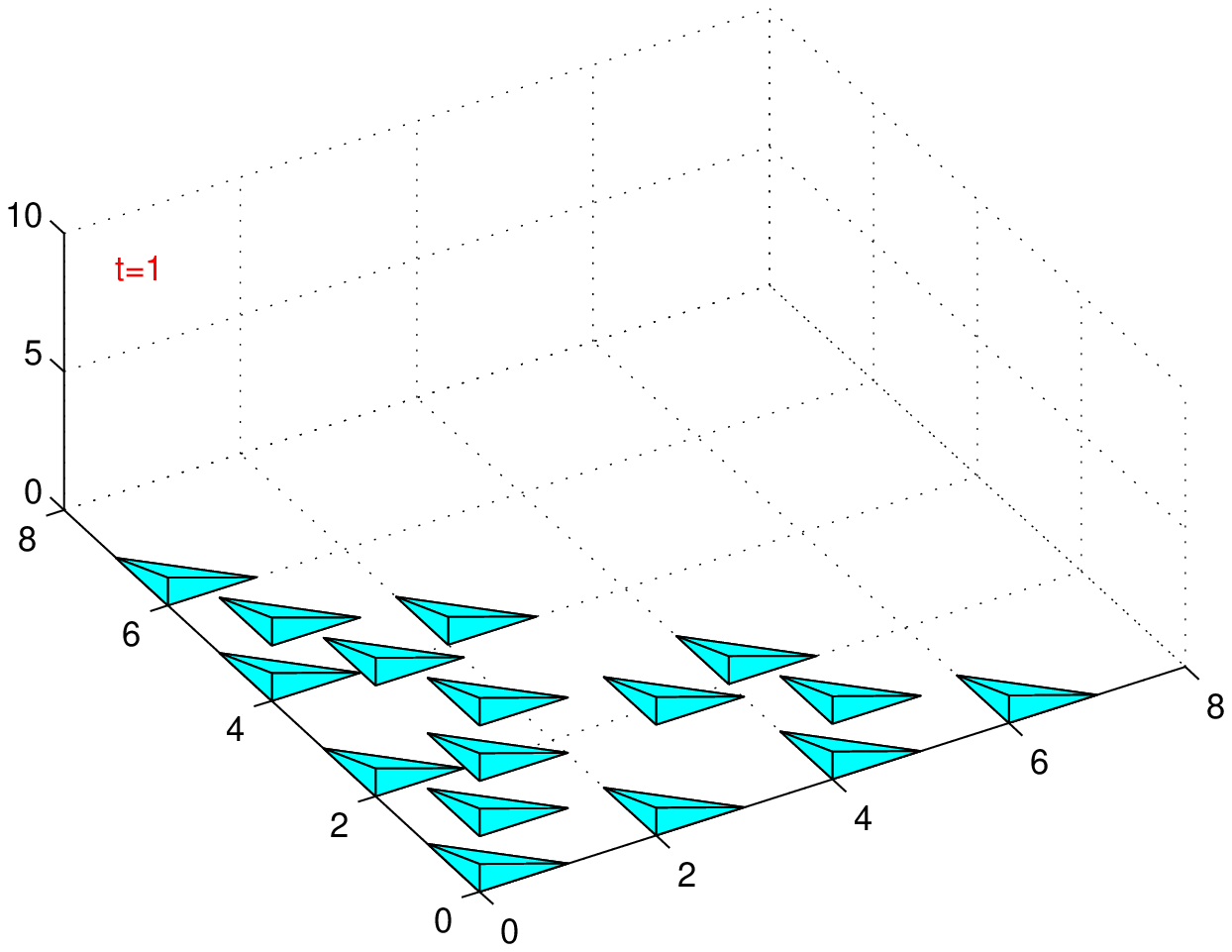}&\includegraphics[width=.3\textwidth]{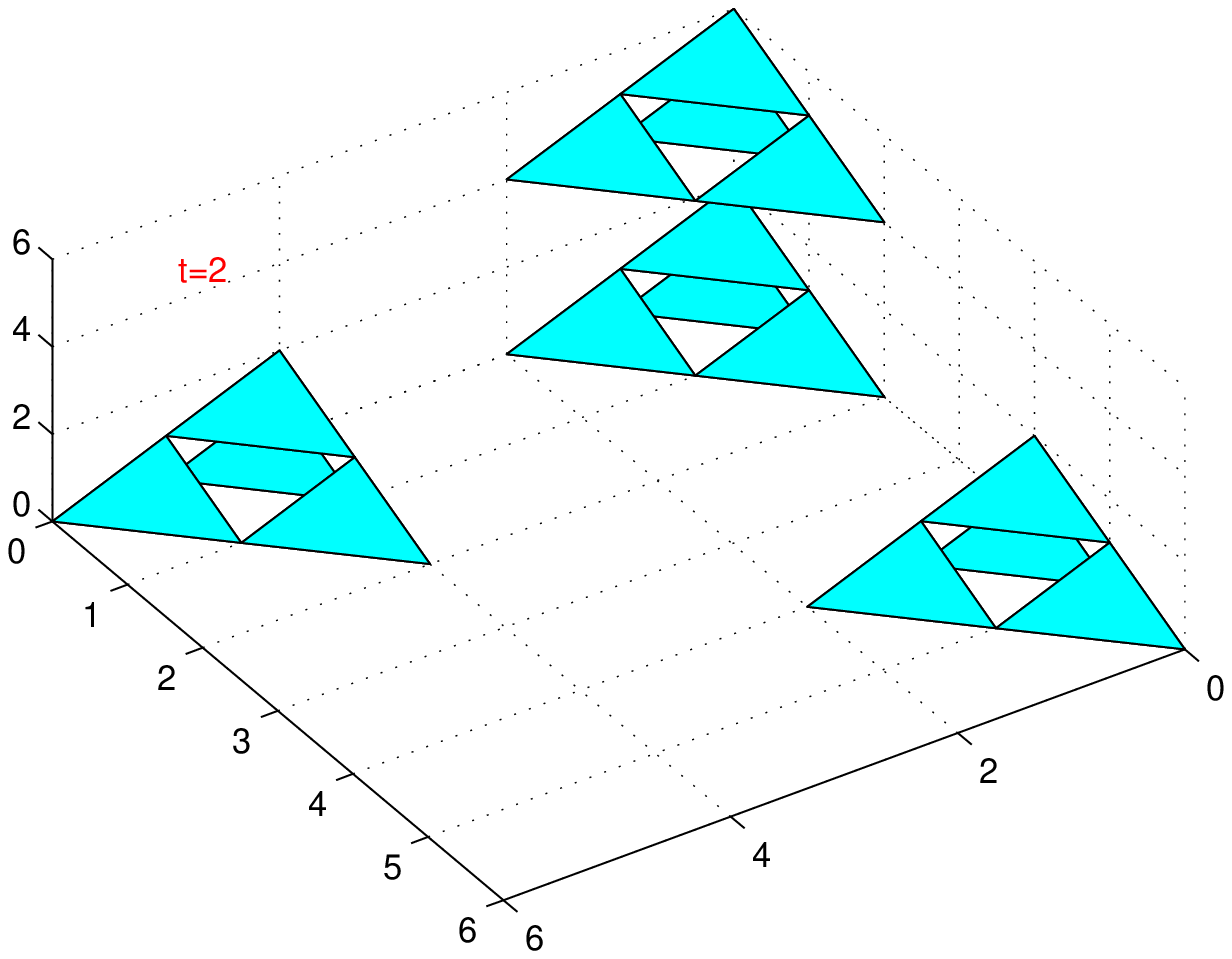}\\
           \includegraphics[width=.3\textwidth]{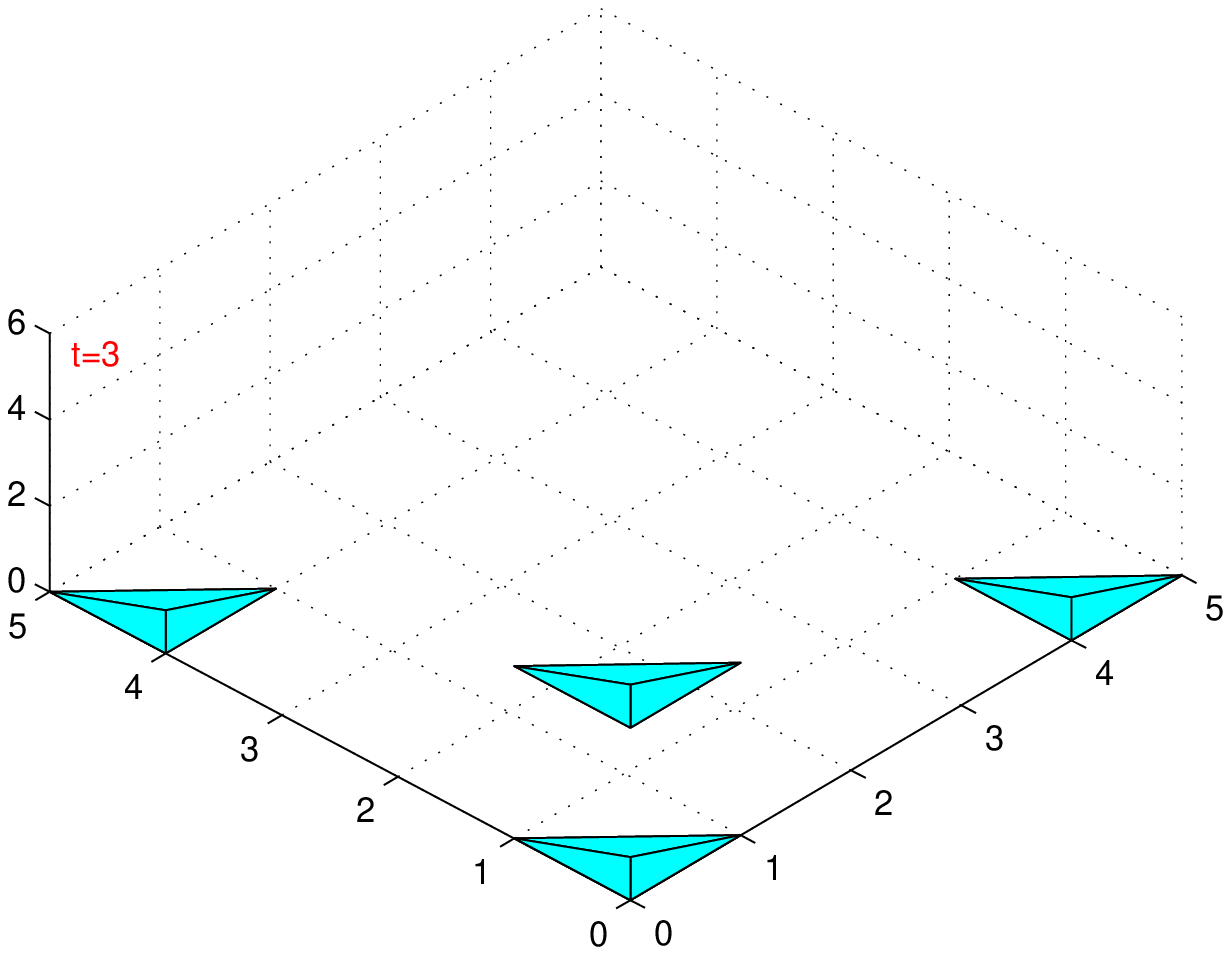}&\includegraphics[width=.3\textwidth]{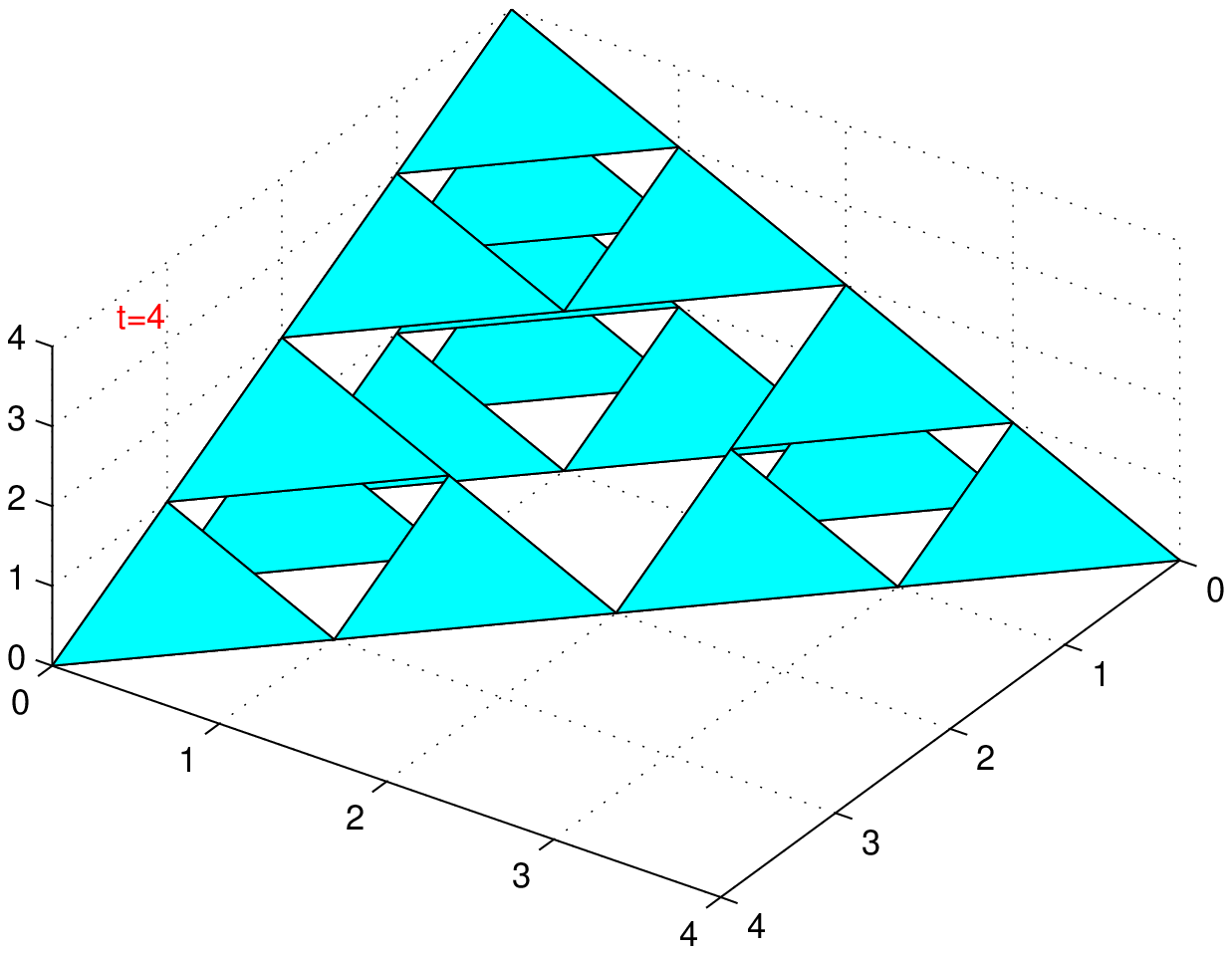}&\includegraphics[width=.3\textwidth]{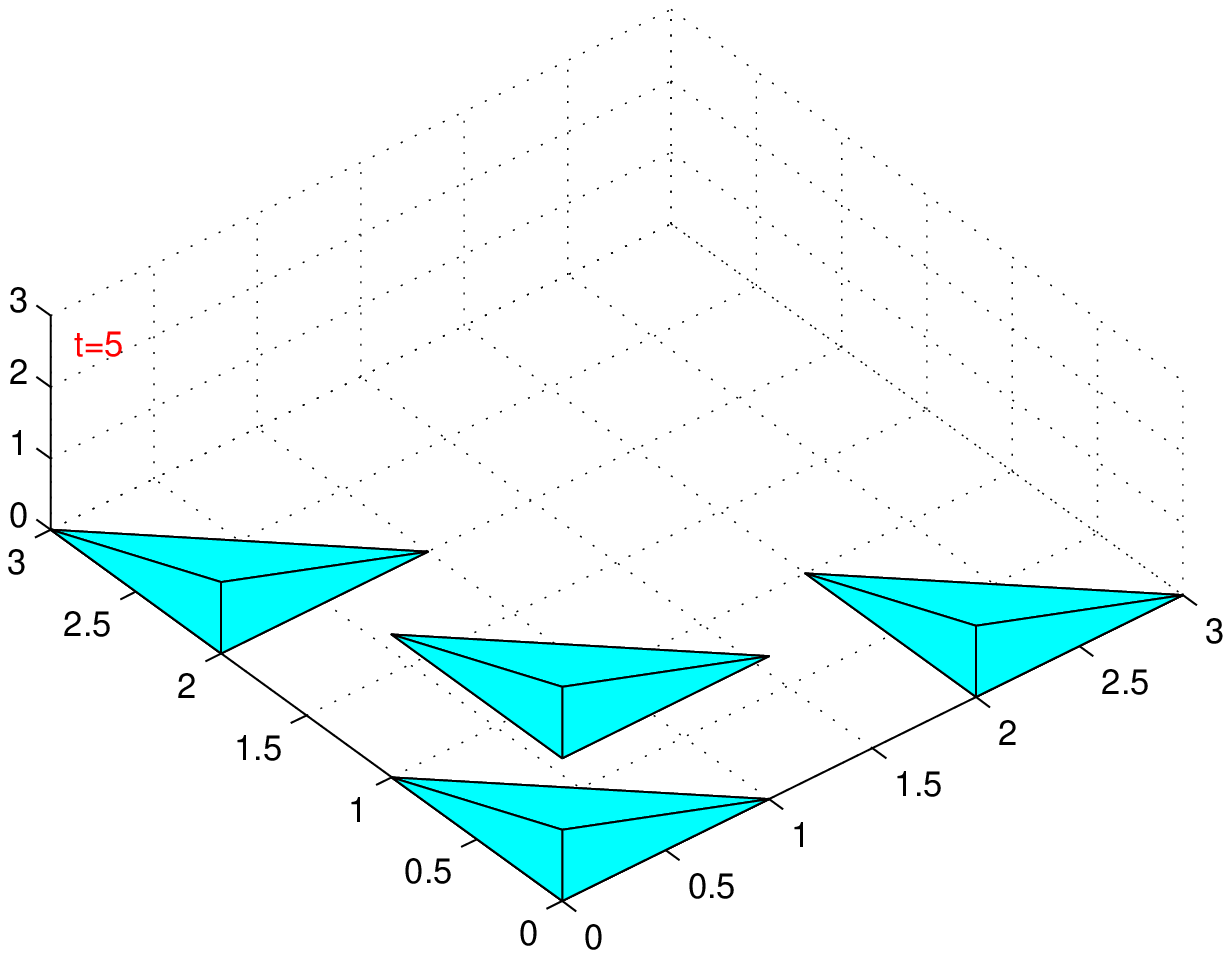}\\
           \includegraphics[width=.3\textwidth]{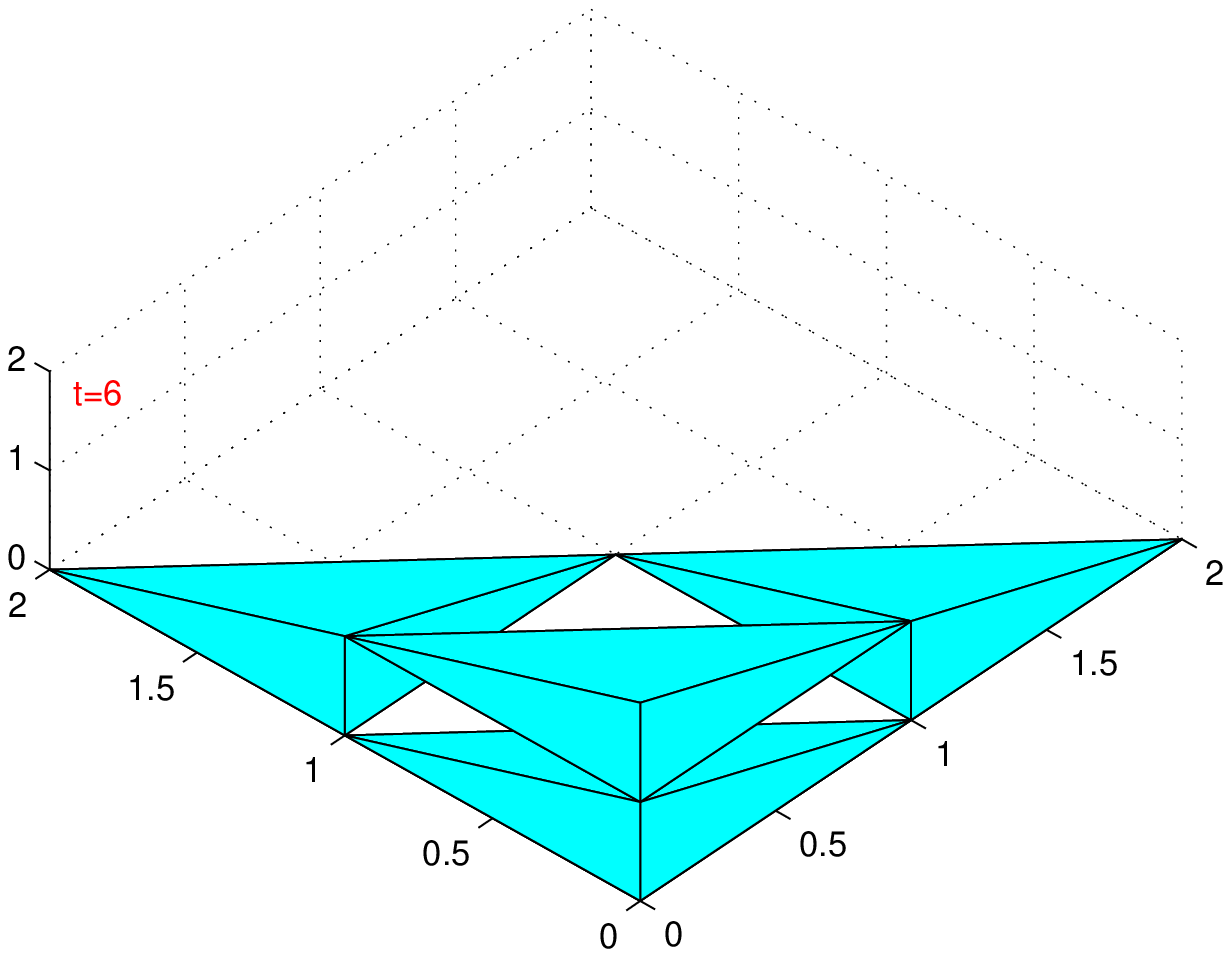}&\includegraphics[width=.3\textwidth]{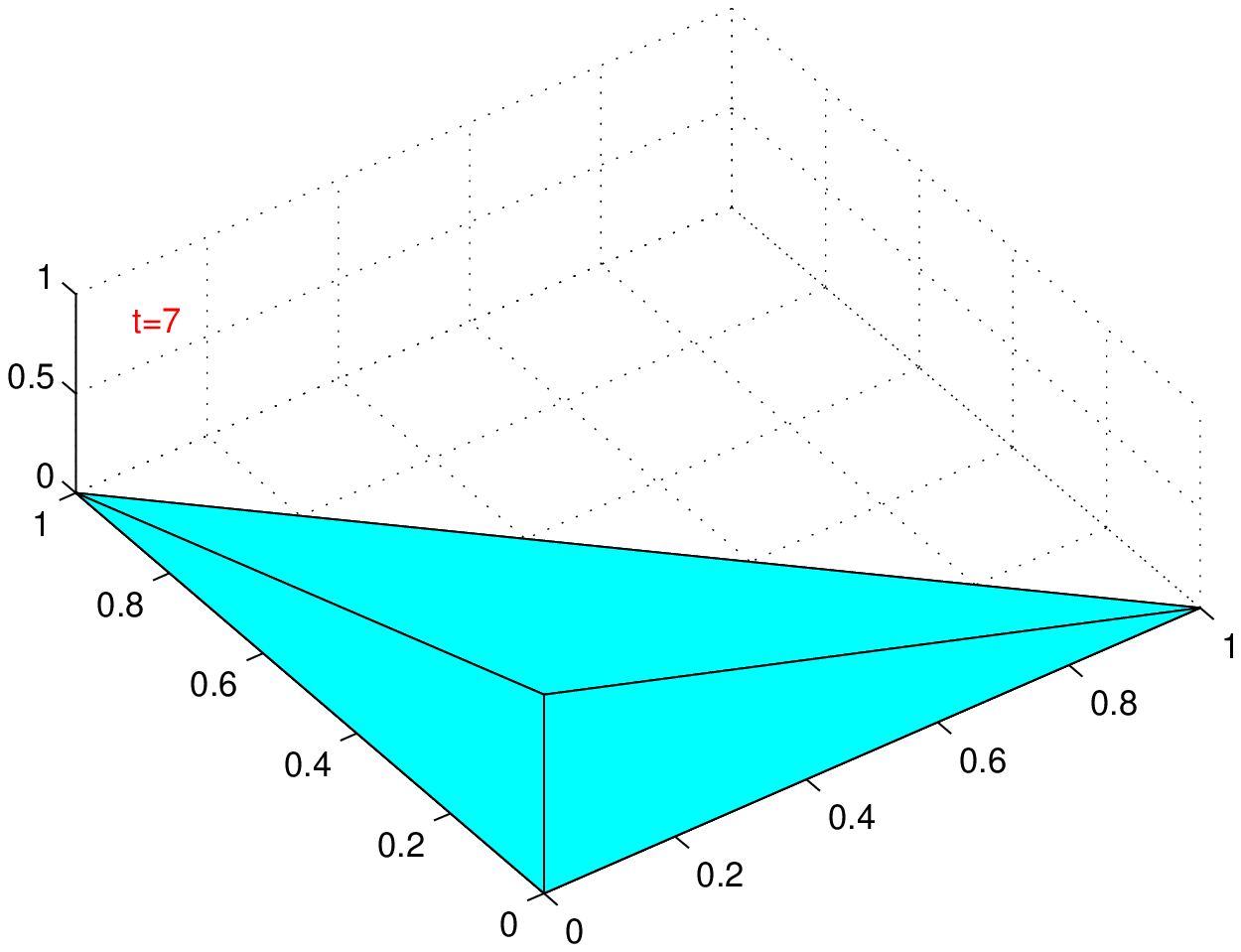}&\includegraphics[width=.3\textwidth]{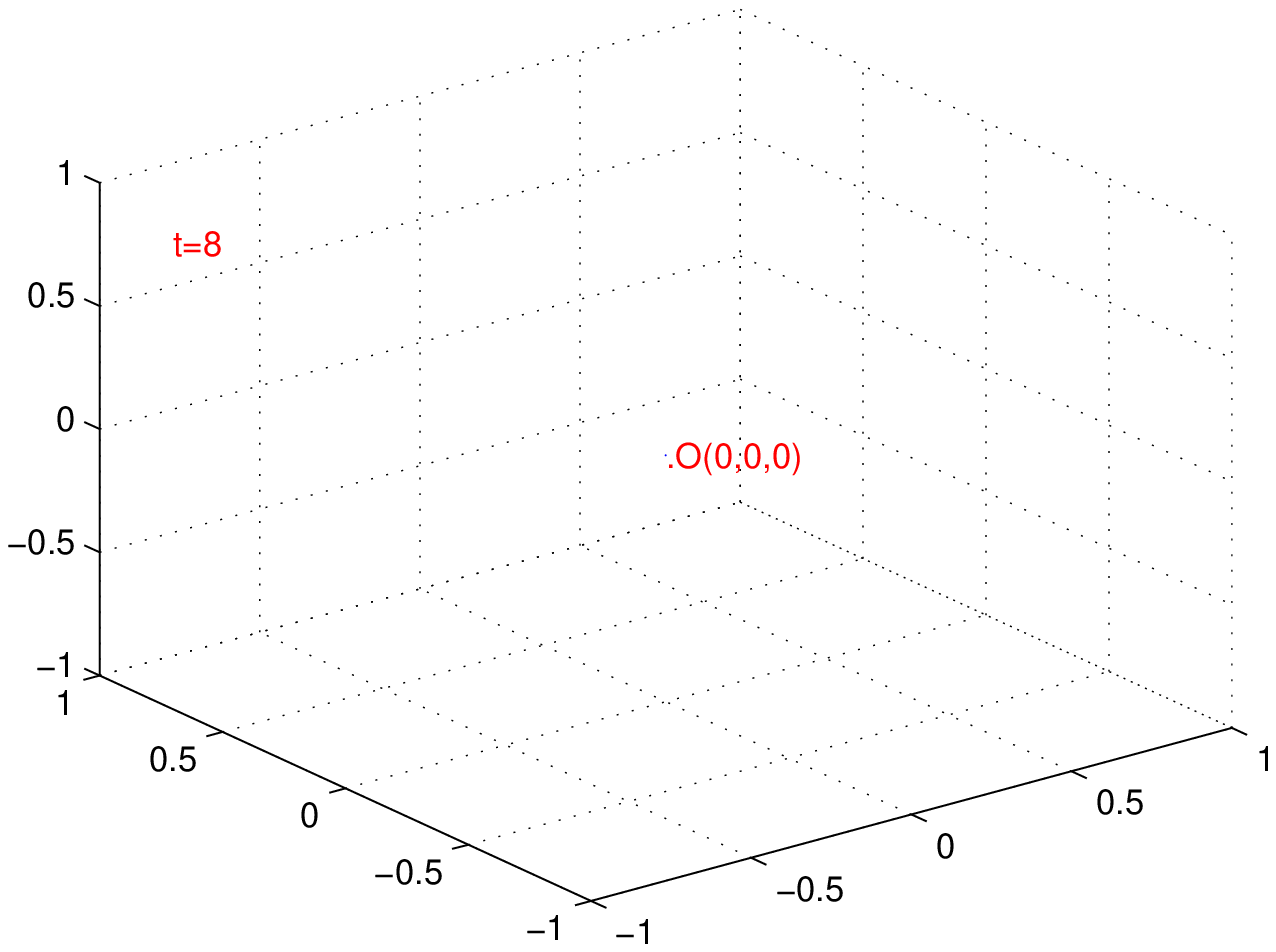}
           \end{tabular}
            \caption{Affine Sierpinski simplex at n=3 in $4$D space with $t=0,1,\cdots,8$. }
            \label{fig-serp-ob3}
 \end{figure}

 Finally, we give the graphs at $n=3$ as shown in Figure \ref{fig-serp-ob3}.
%%%%%%%%%%%%%%%%%%%%%%%%%%%%%%%%
% Begin the references!
%%%%%%%%%%%%%%%%%%%%%%%%%%%%%%%%
\renewcommand{\refname}{\bf\fontsize{12}{12}\selectfont References}
\bibliographystyle{amsplain}

\begin{thebibliography}{99}
\baselineskip 10pt
\bibitem{Bobenko-1} A. Bobenko, T. Hoffmann, B. A. Springborn, \emph{Minimal surfaces from circle patterns:Geometry from combinatorics}, Annals of Mathematics, \textbf{164(1)}(2006), 231-264.
\bibitem{Bobenko-3} A. Bobenko, W. Schief, \emph{Affine spheres:Discreteization via duality relations}, Experimental Mathematics,  \textbf{8(3)}(1999), 261-280.
\bibitem{Bobenko-2} A. Bobenko, P. Schr$\mathrm{\ddot{o}}$der, J. Sullivan, G. Ziegler(Eds.), \emph{Discrete Differential Geometry}, Oberwolfach Seminars, vol. 38, Birkh$\mathrm{\ddot{a}}$user, 2008.
\bibitem{Bobenko-4} A. Bobenko, Y. Suris(Eds.), \emph{Discrete Differential Geometry: Integrable Structure}, Graduate Studies in Mathematics, vol. 98, AMS, 2008.
\bibitem{Brisson} G.F. Brisson, C.A. Reiter, \emph{Sierpinski fractals from words in high dimensions}, Chaos Soliton. Fract. \textbf{5}(1995), 2191-2200.
\bibitem{L-U-Z}A.-M. Li, U. Simon, G. Zhao, \emph{Global affine differential geometry of hypersurfaces}, Berlin - New York: Walter de Gruyter, 1993.
\bibitem{Matsuura} N. Matsuura, H. Urakawa, \emph{Discrete improper affine sphere}, Journal of Geometry and physics,
\textbf{45}(2003), 164-183.
\bibitem{Y-Y-L} Y. Yang, Y. H. Yu, H. L. Liu, \emph{Centroaffine translation surfaces in $\R^3$}, Results in Mathematics, \textbf{56}(2009), 197-210.
\bibitem{Yu-Y-L} Y. H. Yu, Y. Yang, H. L. Liu, \emph{Centroaffine ruled surfaces in $\R^3$}, J. Math. Anal. Appl., \textbf{365}(2010), 683-693.
\bibitem{Wang} C. P. Wang, \emph{Centroaffine minimal hypersurfaces in $\R^{n+1}$}, Geom. Dedicata,
\textbf{51}(1994), 63-74.
\end{thebibliography}

\end{document}